\documentclass[12pt]{article}
\usepackage{amssymb}
\usepackage{amsmath}
\usepackage{color}
\usepackage{comment}
\usepackage{geometry}		% allows easy margin settings
\usepackage{graphicx}
\usepackage[table]{xcolor}  % for coloured tables
\usepackage{tikz}							% for drawing (must be loaded after xcolor)

\let\originalleft\left
\let\originalright\right
\renewcommand{\left}{\mathopen{}\mathclose\bgroup\originalleft}
\renewcommand{\right}{\aftergroup\egroup\originalright}

\newcommand*\circled[1]{\tikz[baseline=(char.base)]{
\node[shape=circle,draw,inner sep=1.2pt](char){#1};}}

\begin{document}

\title{
A Compendium of Hopf-Like Bifurcations in Piecewise-Smooth Dynamical Systems.
}
\author{
D.J.W.~Simpson\\\\
Institute of Fundamental Sciences\\
Massey University\\
Palmerston North\\
New Zealand
}
\maketitle

% keywords:
% MSC codes:

\begin{abstract}

This Letter outlines $20$ geometric mechanisms
by which limit cycles are created locally in two-dimensional piecewise-smooth systems of ODEs.
These include boundary equilibrium bifurcations of hybrid systems, Filippov systems, and continuous systems,
and limit cycles created from folds and by the addition of hysteresis or time-delay.
Scaling laws for the amplitude and period of the limit cycles are compared to (classical) Hopf bifurcations.

\end{abstract}

%Hopf bifurcations are arguably the simplest mechanism by which
%a dynamical system can transition from a steady state to sustained oscillations as parameters varied.
%Hopf bifurcations are perhaps the simplest mechanism by which
%limit cycles can be generated in systems of ODEs.
Hopf bifurcations form perhaps the simplest mechanism
by which limit cycles (isolated periodic orbits) are created in systems of ODEs.
They occur when the real part of a complex conjugate pair of eigenvalues associated with an equilibrium
changes sign as parameters of the system are varied \cite{Ku04}.
%They occur when a complex conjugate pair of eigenvalues associated with an equilibrium
%crosses the imaginary axis in the complex plane as parameters of the system are varied.
%As parameters are varied, a complex conjugate pair of eigenvalues associated with an equilibrium
%cross the imaginary axis in the complex plane.
The limit cycle grows out of the equilibrium
with an amplitude asymptotically proportional to the square root of the parameter change.
The period of the limit cycle varies from $\frac{2 \pi}{\omega}$,
where $\pm {\rm i} \omega$ are eigenvalues of the equilibrium at the bifurcation.
%The period of the limit cycle tends to $\frac{2 \pi}{\omega}$ at the bifurcation,
%where $\pm {\rm i} \omega$ are eigenvalues of the equilibrium.
%Its period is asymptotic to $\frac{2 \pi}{\omega}$, where the aforementioned eigenvalues are $\pm {\rm i} \omega$ at the bifurcation.

In order for Hopf bifurcations to occur in a generic fashion, the ODEs must
be $C^3$ (have continuous third derivatives), at least locally.
Piecewise-smooth ODEs are commonly used to model
physical systems with impacts, switches, or other abrupt processes \cite{DiBu08}.
%For such systems, the lack of smoothness allows limit cycles to be created in different ways.
%For such systems, limit cycles can be created in ways that are not possible for smooth systems.
%For such systems, the presence of codimension-one switching manifolds where the system is nonsmooth
For such systems, the presence of switching manifolds, where the ODEs are not smooth,
allows limit cycles to be created in a wide variety of {\em Hopf-like} bifurcations (HLBs).

%The purpose of this Letter is to briefly summarise and compare different HLBs.
This Letter briefly summarises and compares HLBs.
Details will be provided in a subsequent publication \cite{Si18b}.
For simplicity only two-dimensional systems are treated.
In higher dimensions HLBs are expected to occur in essentially the same way
(with the same scaling laws),
but additional complexities are possible.
%but there is room for additional complexities.
%but certainly more phenomena are possible.

For each type of HLB, one limit cycle is created locally.
Suppose a HLB occurs at $\mu = 0$, where $\mu \in \mathbb{R}$ is a parameter,
and that the limit cycle exists for small $\mu > 0$.
Then there exist constants $k_1, k_2, a > 0$ and $b \ge 0$ such that the amplitude and period obey:
%satisfy the asymptotic formulas
\begin{equation}
\begin{split}
\text{amplitude} &\sim k_1 \mu^a, \\
\text{period} &\sim k_2 \mu^b.
\end{split}
\label{eq:scalingLaws}
\end{equation}
The exponents $a$ and $b$ are determined by the type of HLB bifurcation;
the coefficients $k_1$ and $k_2$ are system specific.
For Hopf bifurcations, $a = \frac{1}{2}$ and $b = 0$.
For a physical system, values for $a$ and $b$ can often be estimated from experimental data.
%From an applied viewpoint, these scaling laws can often be detected from experimental data.
%The results here may aid model selection in that models giving bifurcations
%with incorrect scaling laws could be eliminated.
The results here could aid model selection %the selection of a mathematical model
in that models giving HLBs with incorrect scaling laws would be eliminated.

\begin{figure}[!b]
\begin{center}
\setlength{\unitlength}{1cm}
\begin{picture}(8,3.6)
\put(0,0){\includegraphics[width=3.5cm]{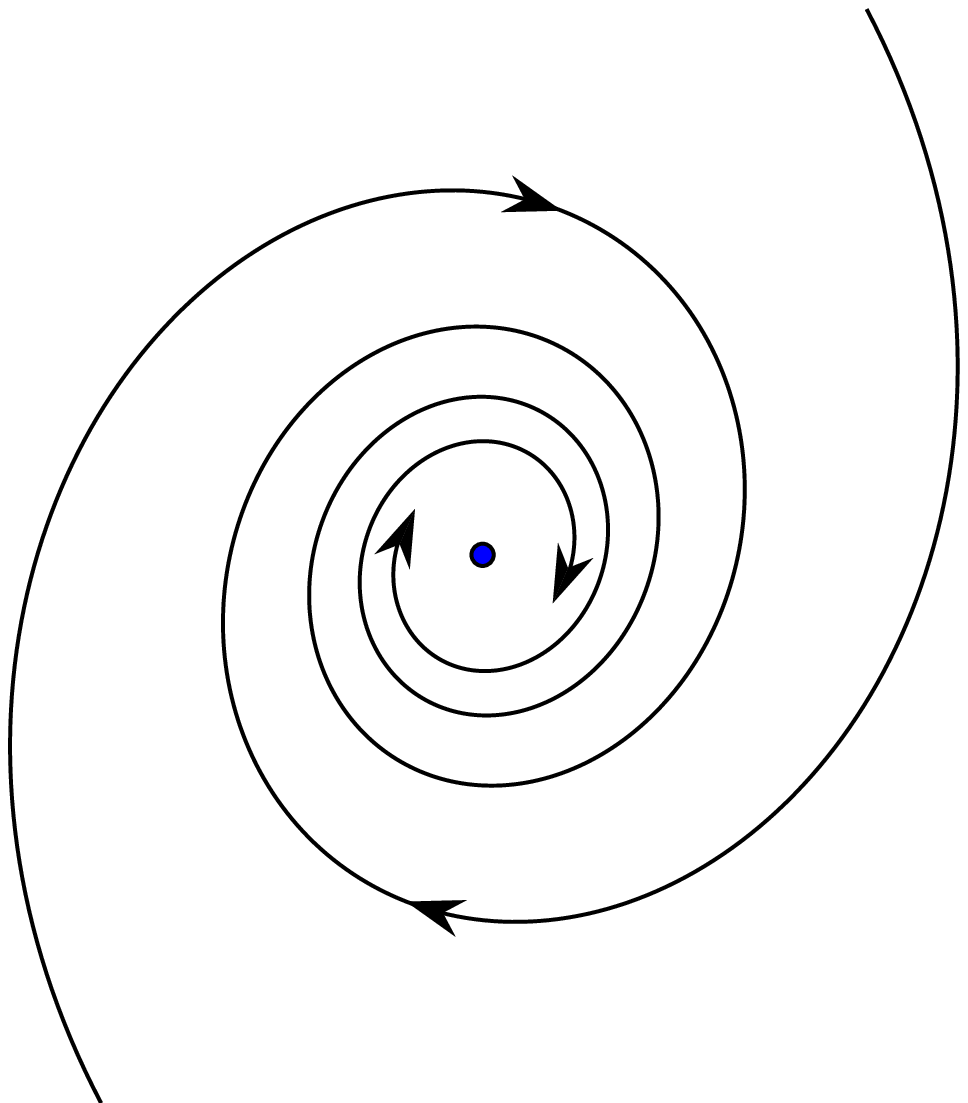}}
\put(4.5,0){\includegraphics[width=3.5cm]{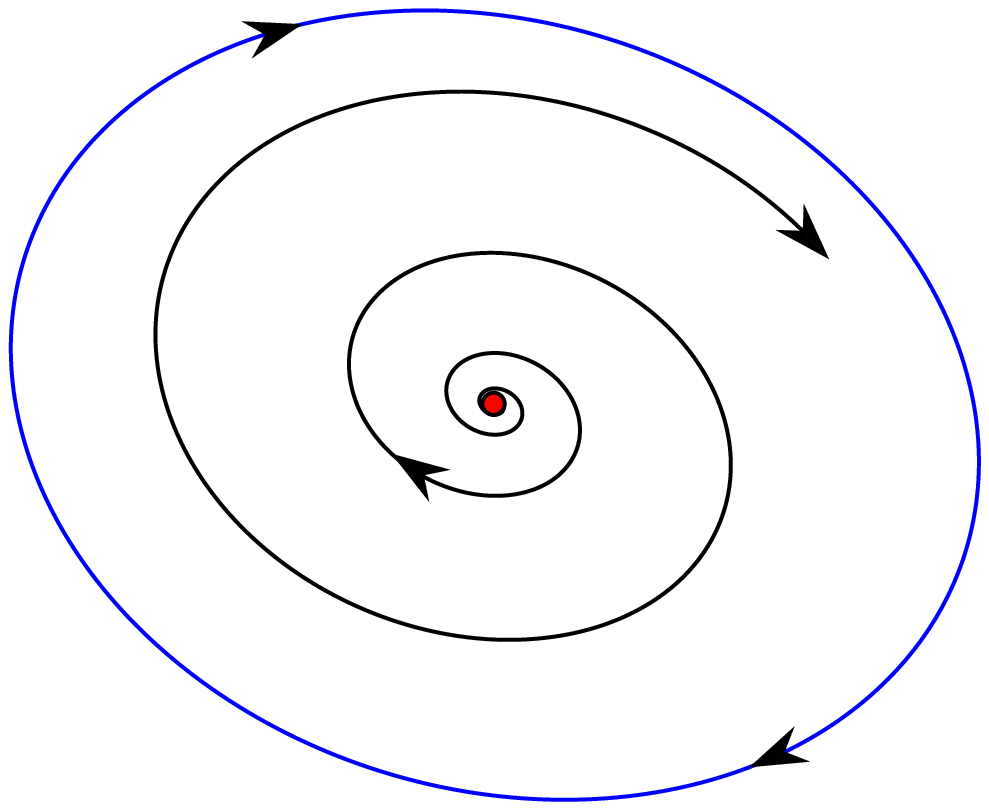}}
\put(0,3.3){\small \parbox{80mm}{\begin{center} $\circled{{\rm H}}$ \end{center}}}
\end{picture}
\caption{
%Phase portraits illustrating a Hopf bifurcation.
Phase portraits of a smooth two-dimensional ODE system
at parameter values either side of a supercritical Hopf bifurcation.
% at which a stable limit cycle is created.
\label{fig:hbl0}
} 
\end{center}
\end{figure}

Table \ref{tb:bifs} lists the HLBs and their values of $a$ and $b$.
%Their numeration in the left-most column is used throughout this Letter.
Mostly $a = 1$ (linear growth) because many of the HLBs are governed by piecewise-linear ODEs.
Indeed, $a \ne 1$ only for HLBs that involve two folds
(a fold is a point on a switching manifold where
one smooth component of the ODEs is tangent to the switching manifold).
%(orbits tangent to a switching manifold).
HLBs with the same values of $a$ and $b$ can be distinguished by qualitative features,
such as the shape of the limit cycle in relation to the switching manifold.
Below, for each type of HLB (numbered 1--20), we give a description and two typical phase portraits
(one for each side of the bifurcation) in cases for which the limit cycle is stable.
Fig.~\ref{fig:hbl0} shows such phase portraits for the Hopf bifurcation.

%%%%%%%%%%%%%%%%%%%%%%%%%%%%%%%%%%%%%%%%%%%%%%%%%%%%%%%%%%%%%
\begin{table}[!t]
\begin{center}
\rowcolors{2}{gray!15}{white}
\setlength\arrayrulewidth{.3mm}
\begin{tabular}{|c|c|c|c|}
\noalign{\hrule height .3mm}
\rowcolor{gray!15}
\# & Description & $a$ & $b$ \\
\noalign{\hrule height .3mm}
H & Hopf & $1/2$ & $0$ \\
$1$ & focus/focus BEB & $1$ & $0$ \\
$2$ & focus/node BEB & $1$ & $0$ \\
$3$ & generic BEB & $1$ & $0$ \\
$4$ & degenerate BEB & $1$ & $0$ \\
$5$ & slipping foci & $1$ & $0$ \\
$6$ & slipping focus/fold & $1$ & $0$ \\
$7$ & slipping folds & $1/2$ & $1/2$ \\
$8$ & fixed foci & $1$ & $0$ \\
$9$ & fixed focus/fold & $1$ & $0$ \\
$10$ & fixed folds & $1/2$ & $1/2$ \\
$11$ & impacting admissible focus & $1$ & $0$ \\
$12$ & impacting virtual focus & $1$ & $0$ \\
$13$ & impacting virtual node & $1$ & $0$ \\
$14$ & impulsive & $1$ & $0$ \\
$15$ & hysteretic pseudo-equilibrium & $1$ & $1$ \\
$16$ & time-delayed pseudo-equilibrium & $1$ & $1$ \\
$17$ & hysteretic two-fold & $1/3$ & $1/3$ \\
$18$ & time-delayed two-fold & $1/2$ & $1/2$ \\
$19$ & intersecting discontinuity surfaces & $1$ & $1$ \\
$20$ & square-root singularity & $1$ & $0$ \\
\noalign{\hrule height .3mm}
\end{tabular}
\caption{
The exponents in the scaling laws \eqref{eq:scalingLaws}
for Hopf bifurcations and $20$ Hopf-like bifurcations.
%Scaling laws \eqref{eq:scalingLaws} for Hopf and Hopf-like bifurcations.
%The exponents in the scaling laws for the amplitude and period
%of the limit cycles created in Hopf bifurcations and $20$ Hopf-like bifurcations.
%Each Hopf-like bifurcation is numbered by its theorem given below.
%Note, BEB abbreviates boundary equilibrium bifurcation.
\label{tb:bifs}
}
\end{center}
\end{table}
%%%%%%%%%%%%%%%%%%%%%%%%%%%%%%%%%%%%%%%%%%%%%%%%%%%%%%%%%%%%%

%%%%%%%%%%%%%%%%%%%%%%%%%%%%%%%%%%%%%%%%%%%%%%%%%%%%%%%%%%%%%%%%%%%%%%%%%%%%%%%%%%%%%%%%%%%%%%%%%%%%%%%%%%%%%%%%%%%%%%%%
%%%%%%%%%%%%%%%%%%%%%%%%%%%%%%%%%%%%%%%%%%%%%%%%%%%%%%%%%%%%%%%%%%%%%%%%%%%%%%%%%%%%%%%%%%%%%%%%%%%%%%%%%%%%%%%%%%%%%%%%

%%%%%%%%%%%%%%%%%%%%%%%%%%%%%%%%%%%%%%%%%%%%%%%%%%%%%%%%%%%%%%%%%%%%%%%%%%%%%%%%%%%%%%%%%%%%%%%%%%%%%%%%%%%%%%%%%%%%%%%%
\begin{figure}[!b]
\begin{center}
\setlength{\unitlength}{1cm}
\begin{picture}(8,15.6)
\put(0,12){\includegraphics[width=3.5cm]{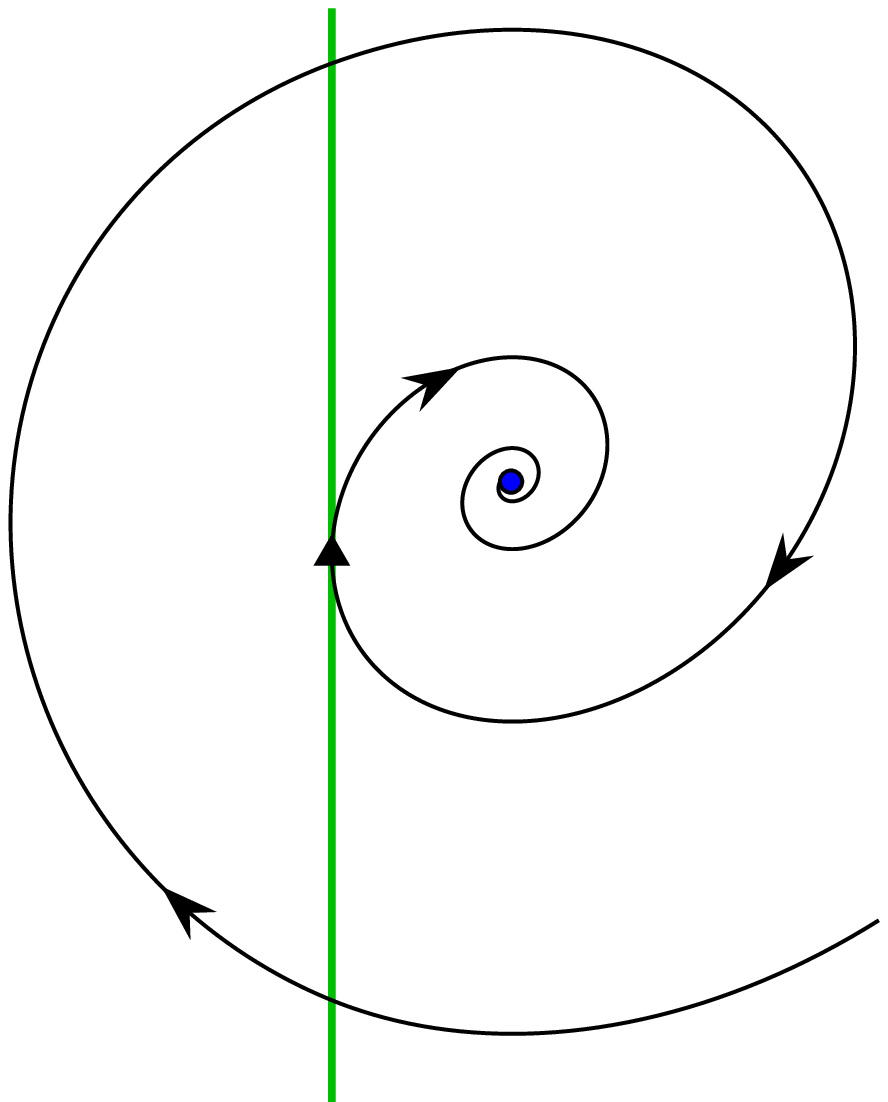}}
\put(4.5,12){\includegraphics[width=3.5cm]{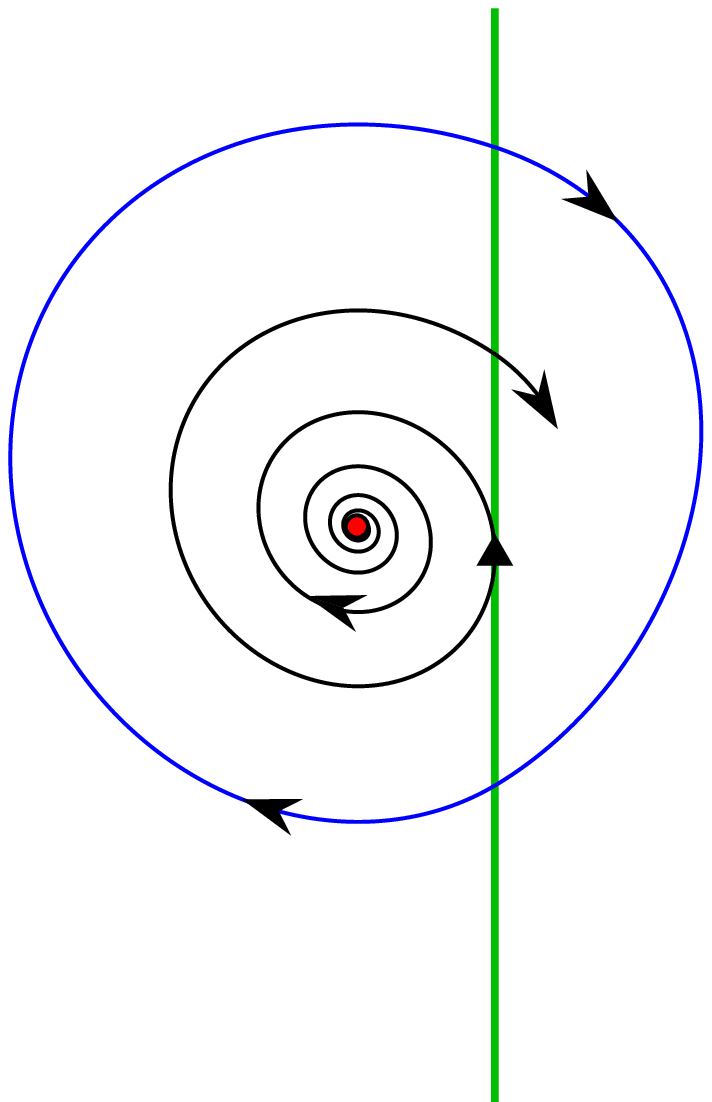}}
\put(0,8){\includegraphics[width=3.5cm]{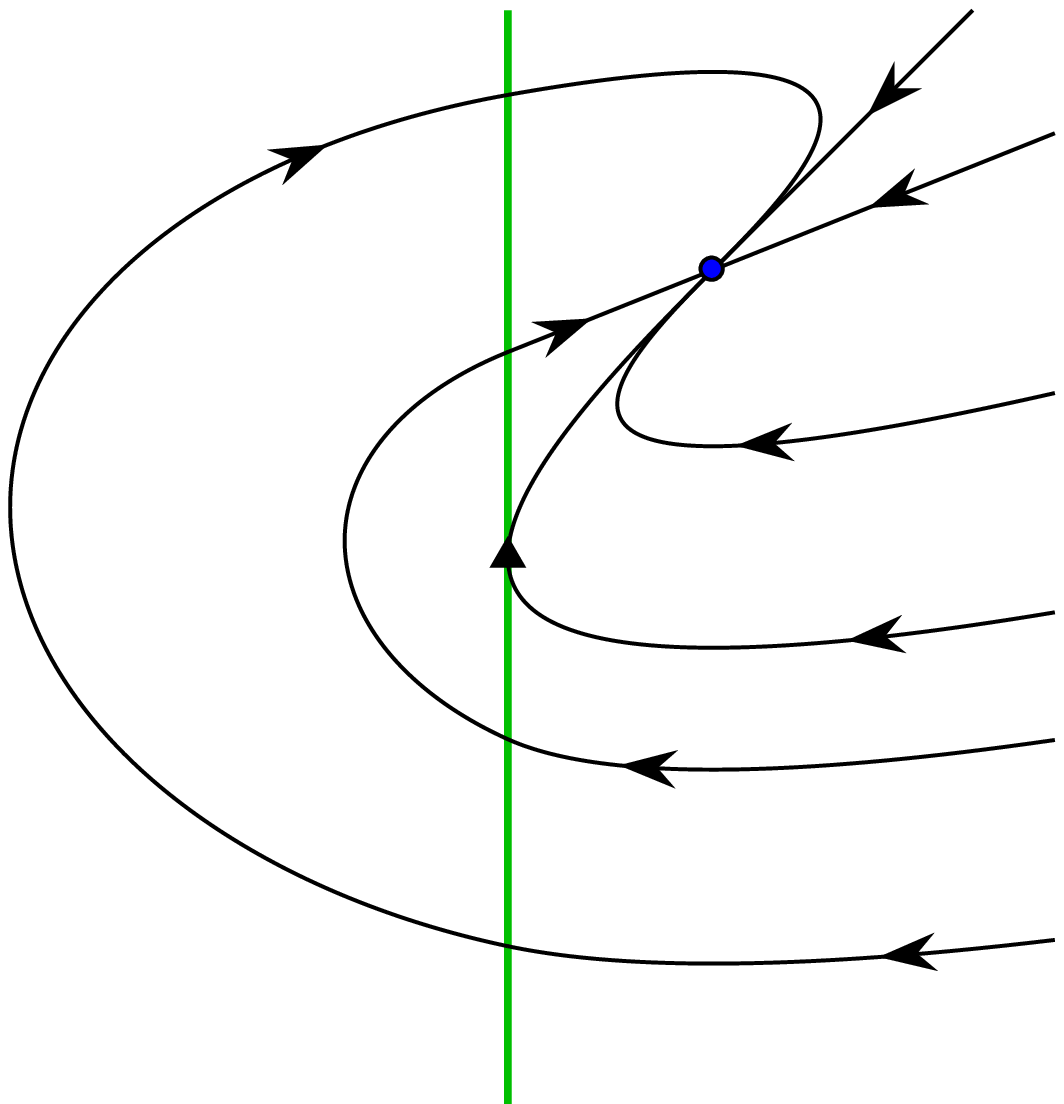}}
\put(4.5,8){\includegraphics[width=3.5cm]{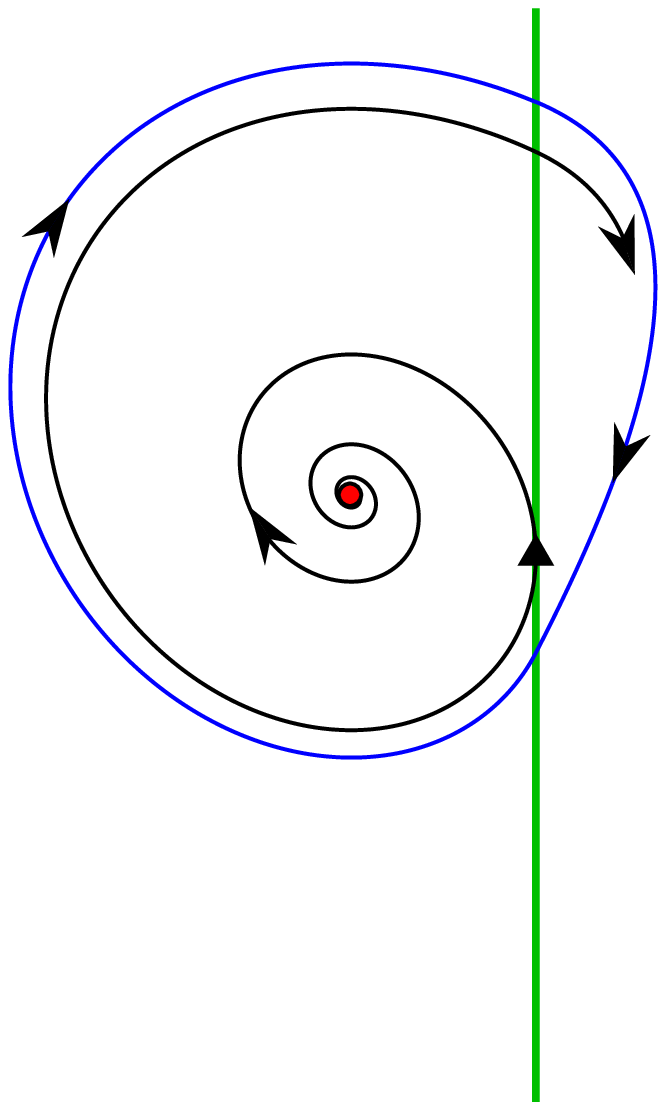}}
\put(0,4){\includegraphics[width=3.5cm]{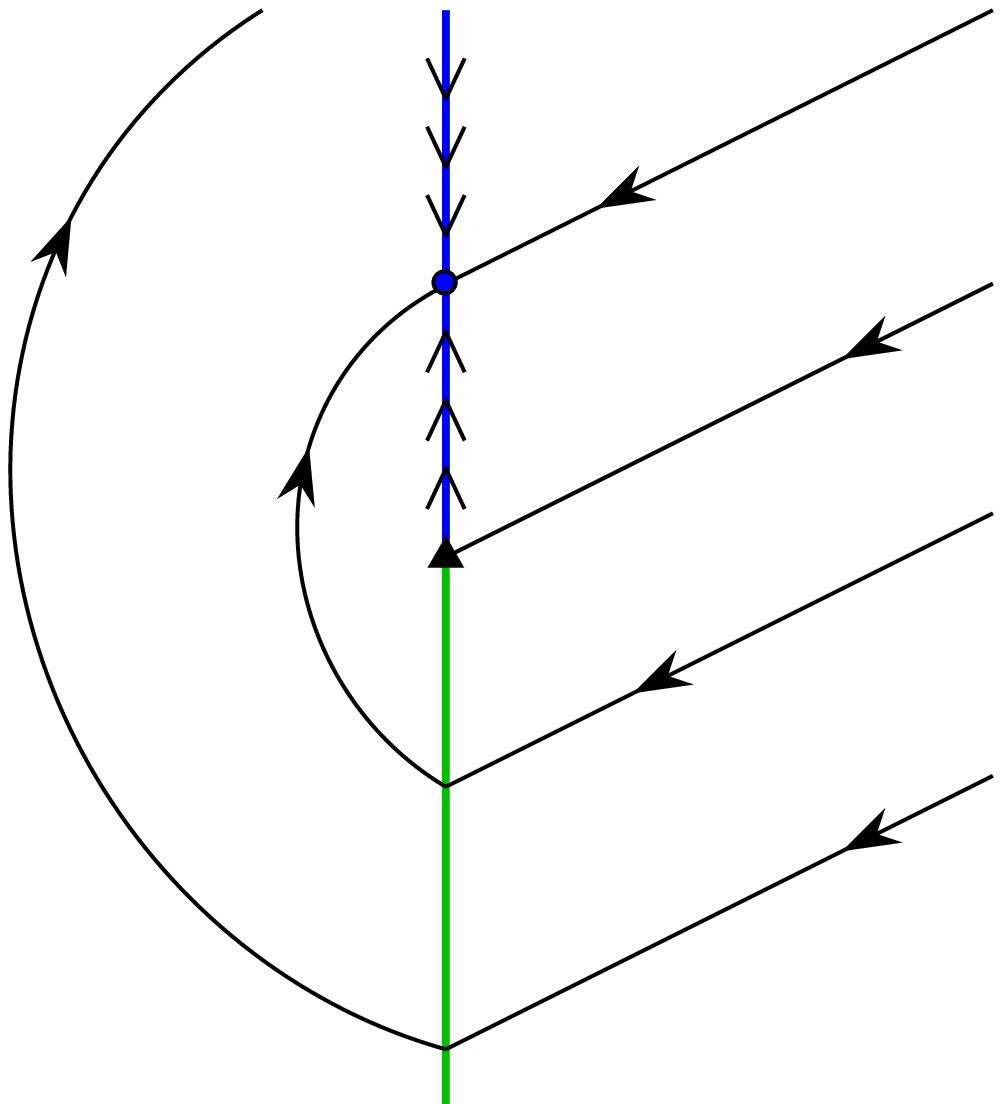}}
\put(4.5,4){\includegraphics[width=3.5cm]{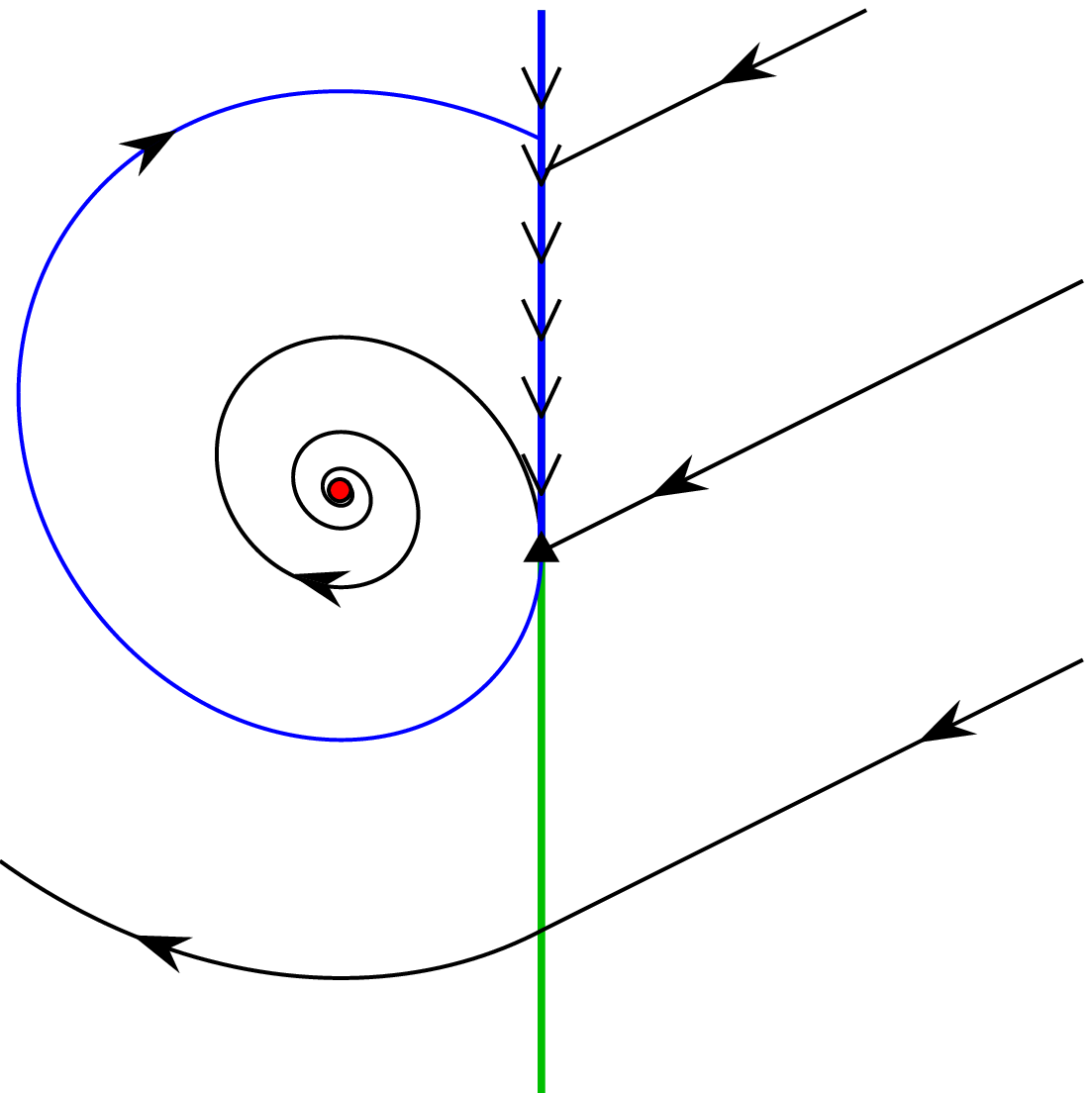}}
\put(0,0){\includegraphics[width=3.5cm]{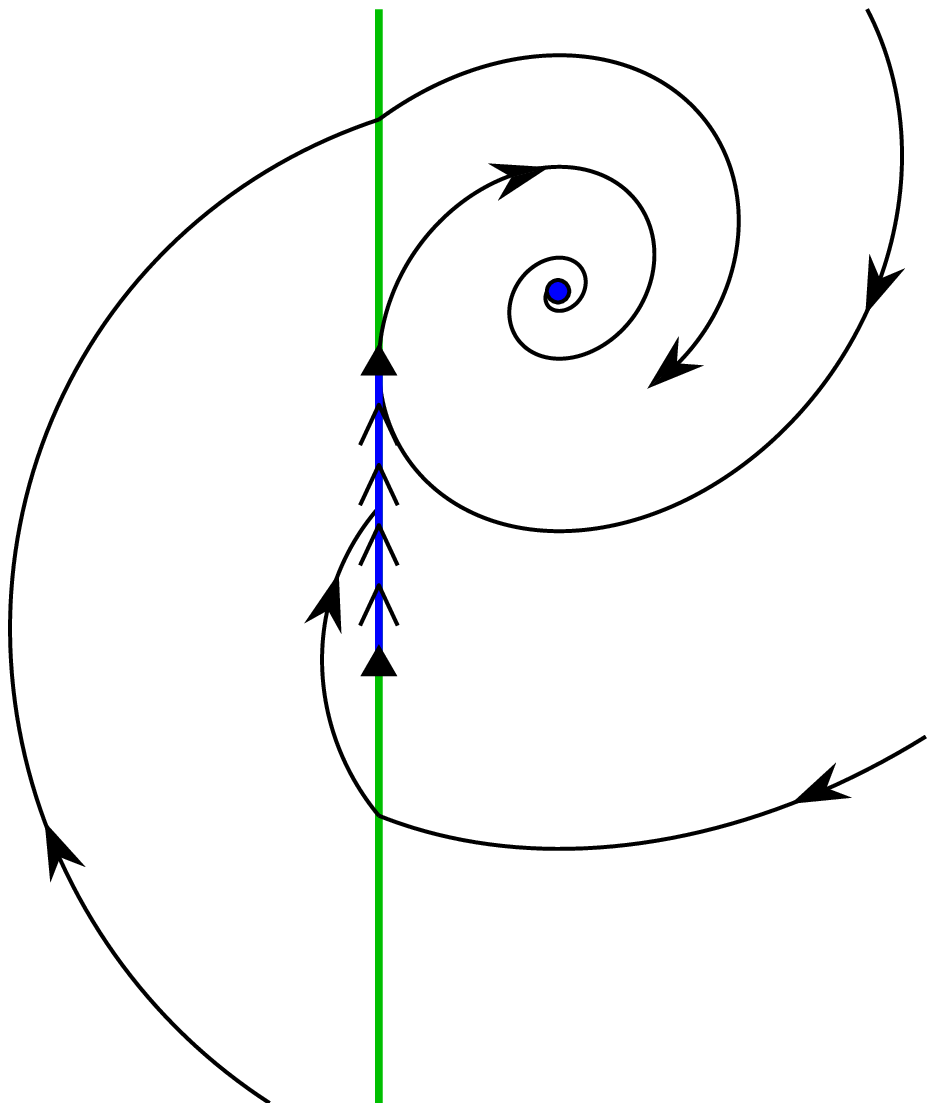}}
\put(4.5,0){\includegraphics[width=3.5cm]{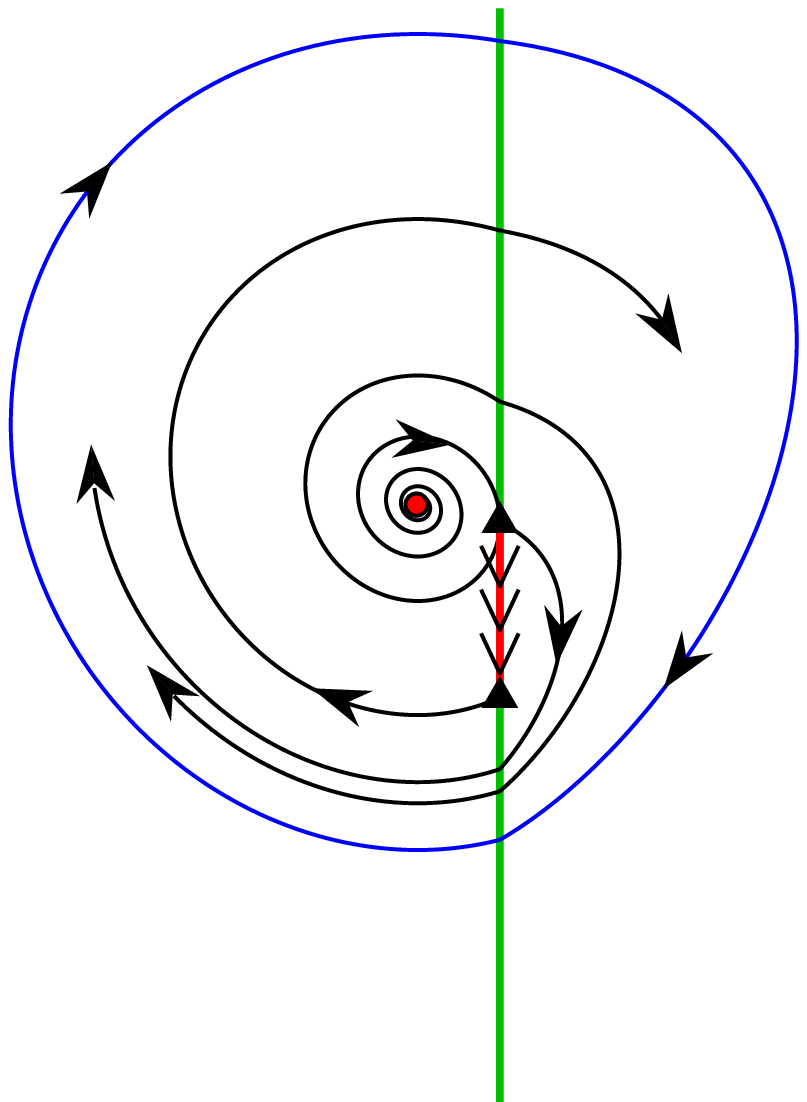}}
\put(0,15.3){\small \parbox{80mm}{\begin{center} $\circled{1}$ \end{center}}}
\put(0,11.3){\small \parbox{80mm}{\begin{center} $\circled{2}$ \end{center}}}
\put(0,7.3){\small \parbox{80mm}{\begin{center} $\circled{3}$ \end{center}}}
\put(0,3.3){\small \parbox{80mm}{\begin{center} $\circled{4}$ \end{center}}}
\end{picture}
\caption{
Phase portraits for boundary equilibrium bifurcations (BEBs).
The switching manifold is indicated by a green line.
Equilibria are shown as circles; folds are shown as triangles.
Stable equilibria, stable limit cycles, and attracting sliding regions are coloured blue.
Unstable equilibria and repelling sliding regions are coloured red.
%Stable [unstable] equilibria are coloured blue [red].
%Similarly, attracting [repelling] sliding regions are coloured blue [red].
\label{fig:hbl1to4}
} 
\end{center}
\end{figure}
%%%%%%%%%%%%%%%%%%%%%%%%%%%%%%%%%%%%%%%%%%%%%%%%%%%%%%%%%%%%%%%%%%%%%%%%%%%%%%%%%%%%%%%%%%%%%%%%%%%%%%%%%%%%%%%%%%%%%%%%

The first four HLBs in Table \ref{tb:bifs} are boundary equilibrium bifurcations (BEBs)
where an equilibrium collides with a switching manifold.
In each case, if the limit cycle is stable, it must encircle an unstable focus.
%First suppose the system is continuous on the switching manifold and, at least locally,
First, consider a system that is continuous on a switching manifold, and, at least locally,
can be put in the form
\begin{equation}
\begin{bmatrix} \dot{x} \\ \dot{y} \end{bmatrix} =
\begin{cases} F_L(x,y;\mu), & x \le 0, \\ F_R(x,y;\mu), & x \ge 0, \end{cases}
\label{eq:pwscODE}
\end{equation}
where $x$ and $y$ are the state variables and $x=0$ is the switching manifold.
By assumption, \eqref{eq:pwscODE} is continuous but non-differentiable on $x=0$,
so at the BEB the eigenvalues associated with the equilibrium typically change discontinuously.
For HLB 1, see Fig.~\ref{fig:hbl1to4}, the equilibrium changes from an unstable focus to a stable focus.
In order for a stable limit cycle to be created,
the attraction of the stable focus must dominate the repulsion of the unstable focus.
Specifically we need $\alpha < 0$, where
\begin{equation}
\alpha = \frac{\lambda_L}{\omega_L} + \frac{\lambda_R}{\omega_R},
\label{eq:pwscNondegCond}
\end{equation}
and $\lambda_L \pm {\rm i} \omega_L$ and $\lambda_R \pm {\rm i} \omega_R$,
with $\lambda_L > 0$, $\lambda_R < 0$, $\omega_L > 0$, and $\omega_R > 0$,
are the eigenvalues associated with the equilibria at the bifurcation \cite{FrPo97,SiMe07}.
For HLB 2, the equilibrium changes to a stable node
and there is no such criticality condition \cite{FrPo98,SiMe12}.
Both bifurcations are governed by the linear terms
in a piecewise expansion of \eqref{eq:pwscODE}, and so $a = 1$.
Also, $b = 0$, but unlike Hopf bifurcations
the limiting value of the period is not given by a simple expression.
These HLBs have been identified in the McKean neuron model \cite{Mc70}
and other piecewise-linear models of excitable systems \cite{DeFr13,RoCo12}.

Next we consider Filippov systems of the form
\begin{equation}
\begin{bmatrix} \dot{x} \\ \dot{y} \end{bmatrix} =
\begin{cases} F_L(x,y;\mu), & x < 0, \\ F_R(x,y;\mu), & x > 0, \end{cases}
\label{eq:FilippovBEB}
\end{equation}
which are discontinuous on $x=0$.
Subsets of $x=0$ at which $F_L$ and $F_R$ both point towards $x=0$
%[away from $x=0$] are attracting [repelling] sliding regions.
are attracting sliding regions.
When an orbit reaches an attracting sliding region (as time increases),
it subsequently evolves on $x=0$.
For Filippov systems such sliding motion is governed by the convex combination of $F_L$ and $F_R$
tangent to $x=0$ \cite{Fi88}.

When an unstable focus of \eqref{eq:FilippovBEB} collides with $x=0$ in a BEB,
it may turn into an attracting pseudo-equilibrium (an equilibrium of the sliding motion). 
In this case a limit cycle exists with the focus (HLB 3).
The limit cycle involves only one side of the switching manifold and has a segment of sliding motion.
This type of HLB occurs in, for instance, the Gause predator-prey model \cite{Kr11}. % [GaSm36]
%Other BEBs are described in \cite{KuRi03,HoHo16,Gl16d}.

Generic codimension-one BEBs in Filippov systems
involve one equilibrium and one pseudo-equilibrium \cite{KuRi03}. %{HoHo16,Gl16d}.
For symmetric Filippov systems of the form \eqref{eq:FilippovBEB}, the point on $x=0$ at which the bifurcation occurs
can be an equilibrium of both $F_L$ and $F_R$, such as for a circuit system given in \cite{AnVi66}.
Such BEBs resemble that of continuous systems, except sliding motion is possible.
If an unstable focus transitions to a stable focus with $\alpha < 0$,
both foci have the same direction of rotation,
%(where $\alpha$ is given by \eqref{eq:pwscNondegCond}),
and, at least locally, the system has no attracting sliding regions when the unstable focus exists,
then a unique stable limit cycle exists around the unstable focus (HLB 4).
If attracting sliding regions are present,
up to three nested limit cycles may be created at the BEB simultaneously \cite{BrMe13,FrPo14}.

%%%%%%%%%%%%%%%%%%%%%%%%%%%%%%%%%%%%%%%%%%%%%%%%%%%%%%%%%%%%%%%%%%%%%%%%%%%%%%%%%%%%%%%%%%%%%%%%%%%%%%%%%%%%%%%%%%%%%%%%
\begin{figure}[!b]
\begin{center}
\setlength{\unitlength}{1cm}
\begin{picture}(8,11.6)
\put(0,8){\includegraphics[width=3.5cm]{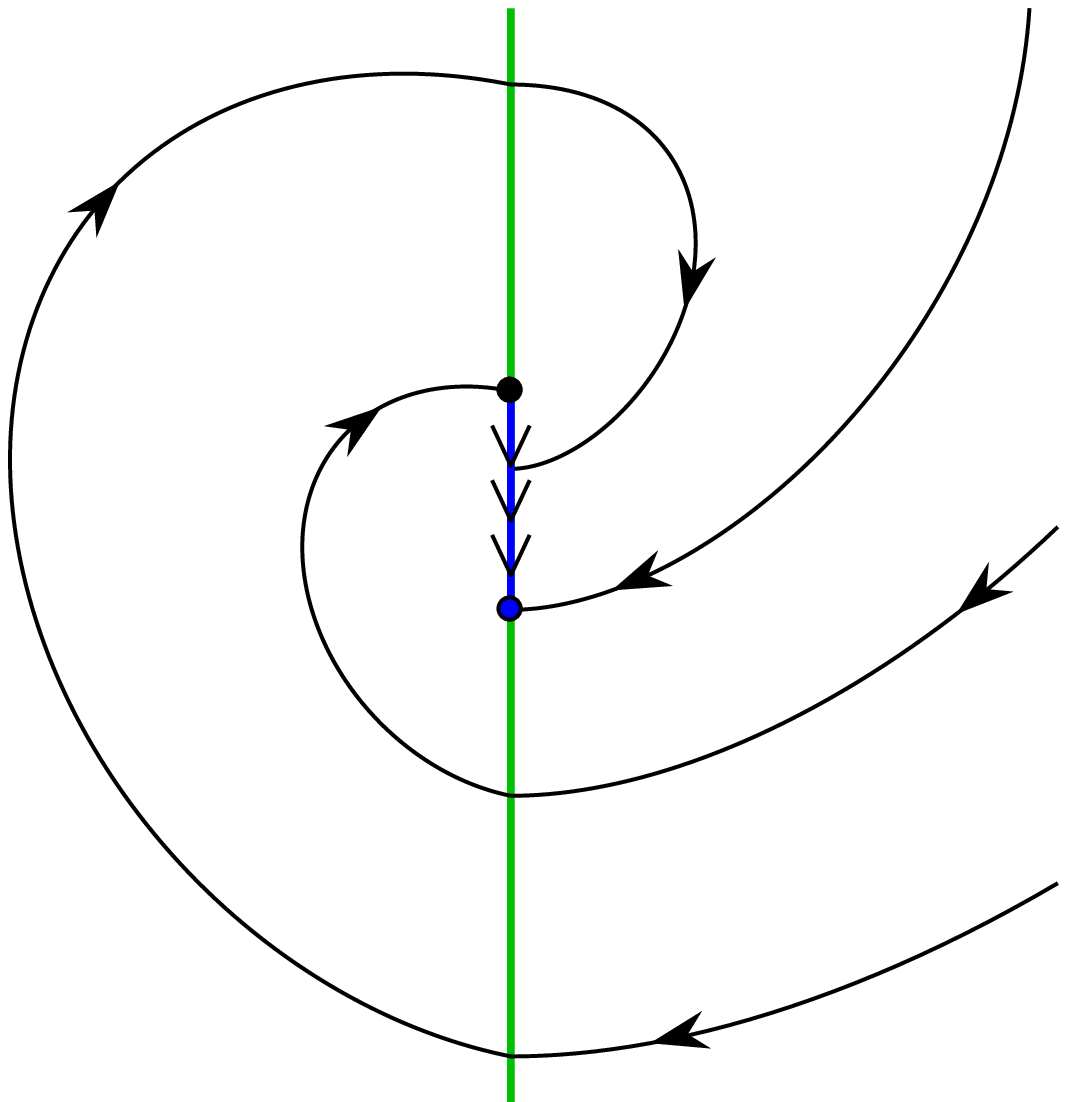}}
\put(4.5,8){\includegraphics[width=3.5cm]{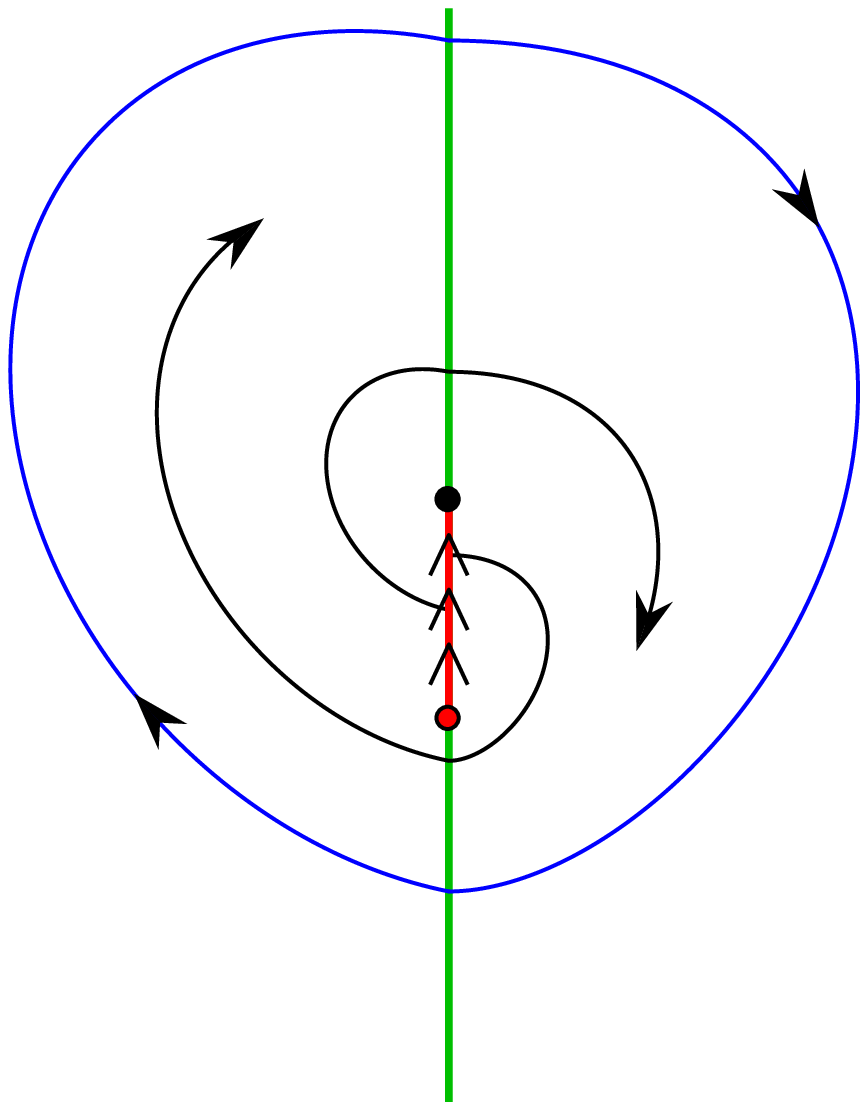}}
\put(0,4){\includegraphics[width=3.5cm]{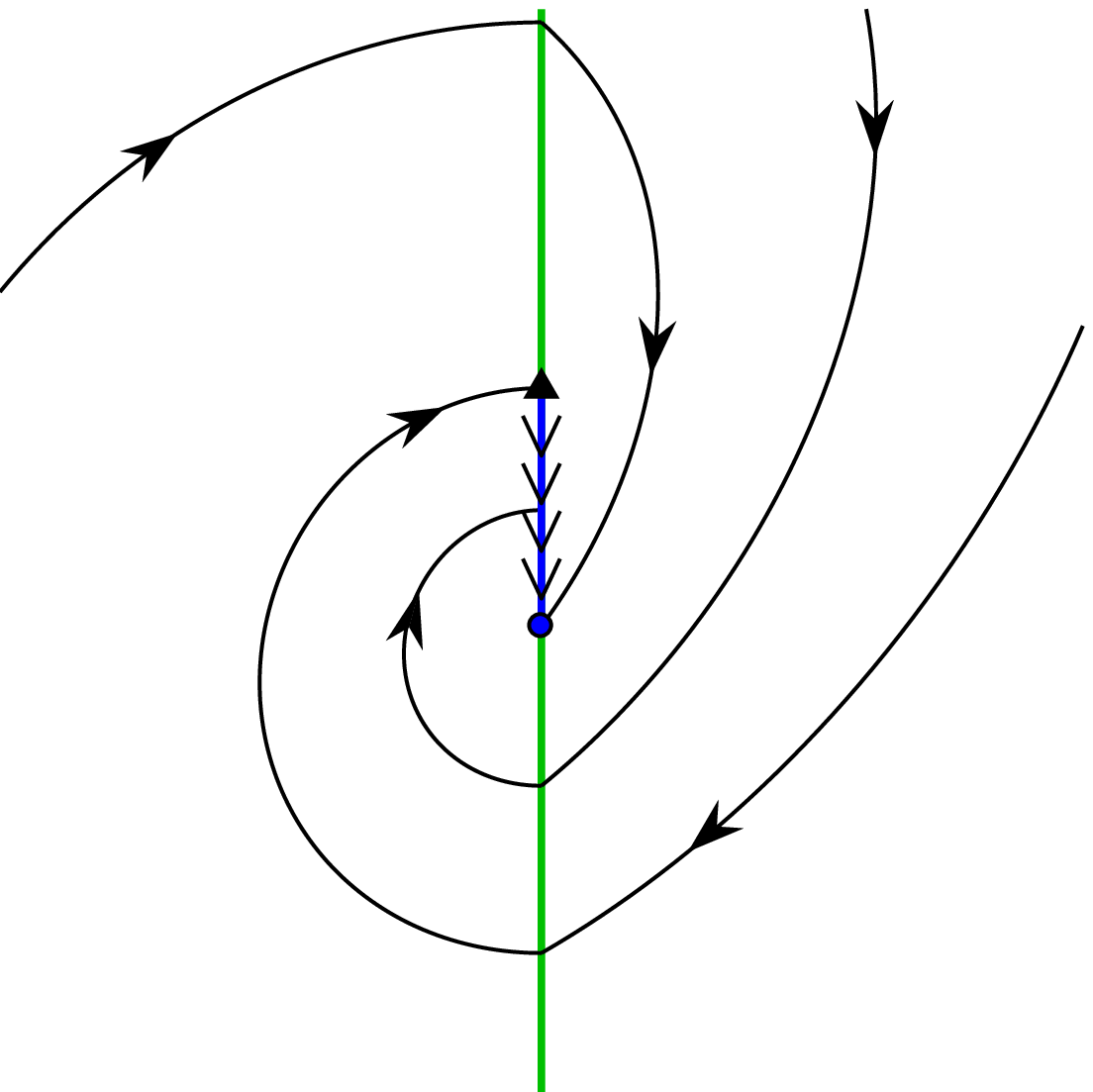}}
\put(4.5,4){\includegraphics[width=3.5cm]{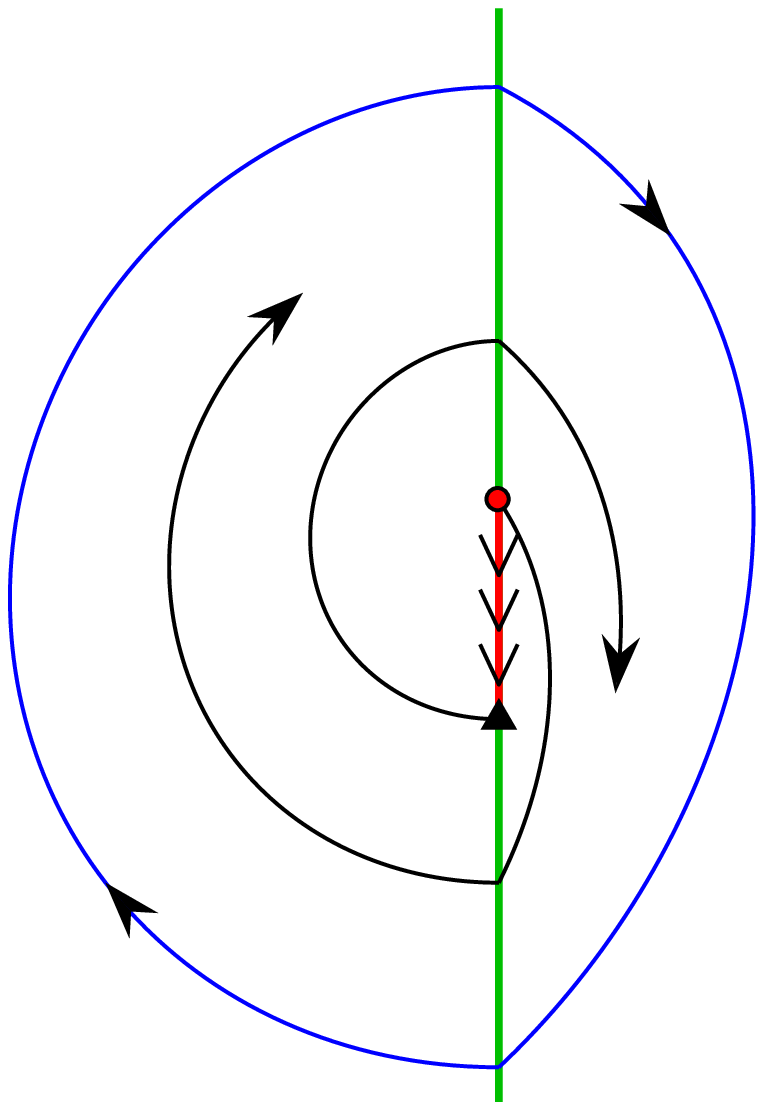}}
\put(0,0){\includegraphics[width=3.5cm]{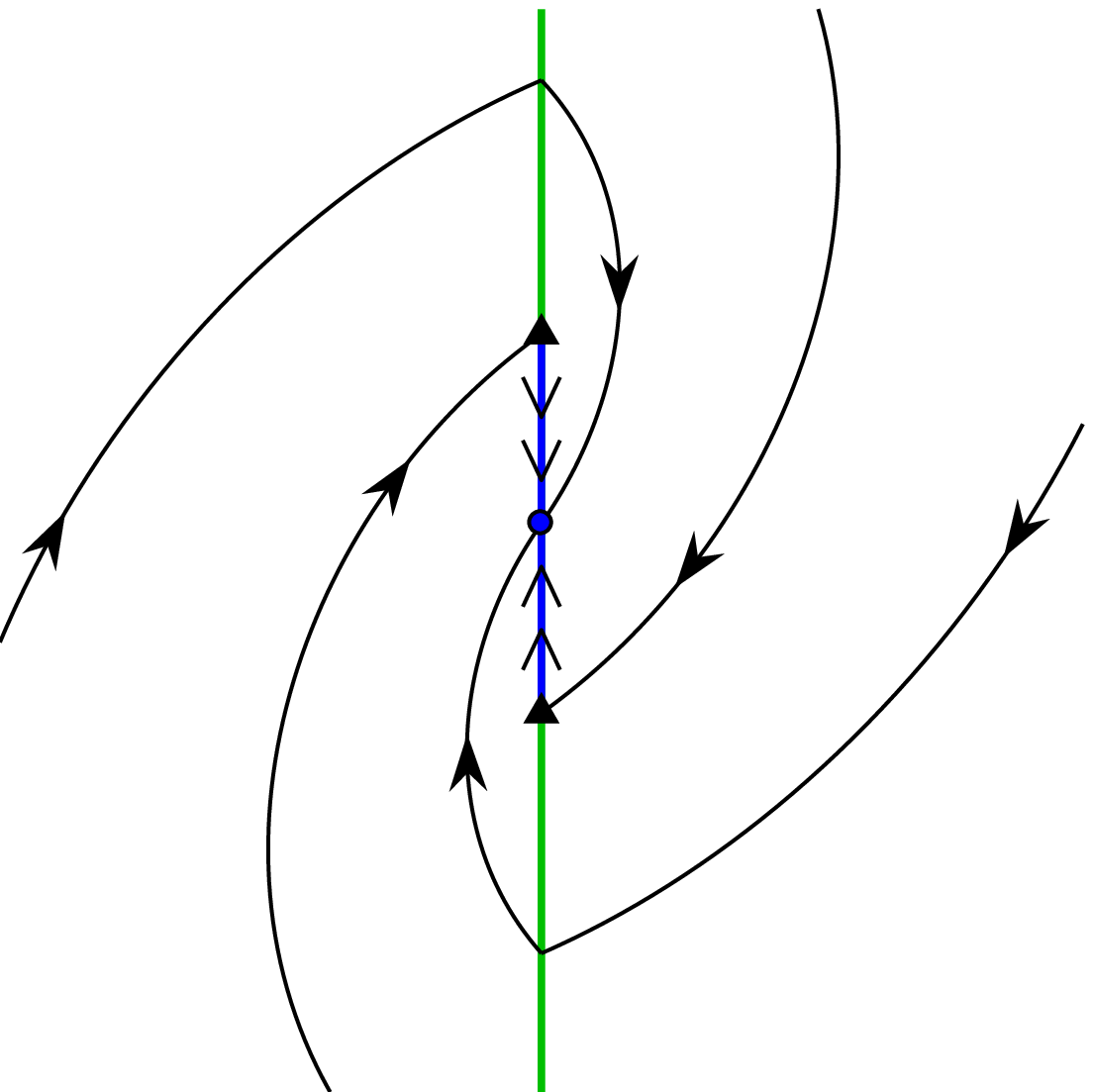}}
\put(4.5,0){\includegraphics[width=3.5cm]{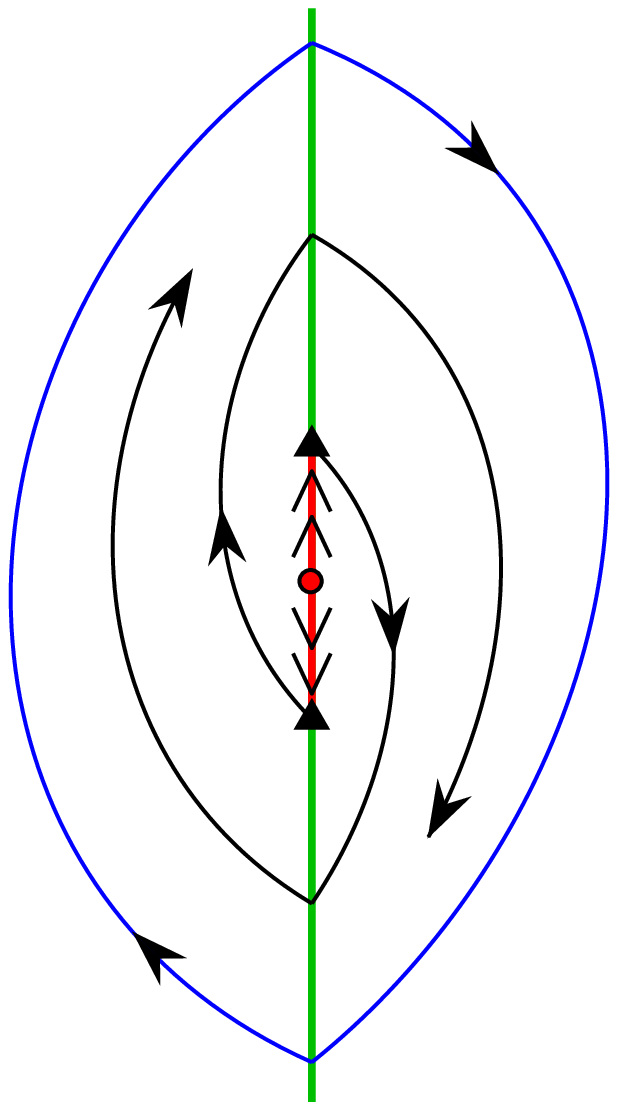}}
\put(0,11.3){\small \parbox{80mm}{\begin{center} $\circled{5}$ \end{center}}}
\put(0,7.3){\small \parbox{80mm}{\begin{center} $\circled{6}$ \end{center}}}
\put(0,3.3){\small \parbox{80mm}{\begin{center} $\circled{7}$ \end{center}}}
\end{picture}
\caption{
Phase portraits for slipping folds and foci.
\label{fig:hbl5to7}
} 
\end{center}
\end{figure}
%%%%%%%%%%%%%%%%%%%%%%%%%%%%%%%%%%%%%%%%%%%%%%%%%%%%%%%%%%%%%%%%%%%%%%%%%%%%%%%%%%%%%%%%%%%%%%%%%%%%%%%%%%%%%%%%%%%%%%%%

Again consider \eqref{eq:FilippovBEB},
but now suppose $F_L$ and $F_R$ each have either a focus or a fold on $x=0$
for all values of $\mu$ in a neighbourhood of $0$.
Furthermore, suppose that the foci/folds `slip' along $x=0$ as the value of $\mu$ is varied,
and collide at $\mu = 0$, Fig.~\ref{fig:hbl5to7}.
A local limit cycle can be created at $\mu = 0$, and there are three cases:
(i) two foci (HLB 5),
(ii) one focus and one fold (HLB 6), and
(iii) two folds (HLB 7).
In the last case the amplitude of the limit cycle is asymptotically proportional to $\sqrt{\mu}$.
This bifurcation is generic and occurs in
a prototypical model of balancing via on-off control \cite{Ko17}.

%%%%%%%%%%%%%%%%%%%%%%%%%%%%%%%%%%%%%%%%%%%%%%%%%%%%%%%%%%%%%%%%%%%%%%%%%%%%%%%%%%%%%%%%%%%%%%%%%%%%%%%%%%%%%%%%%%%%%%%%
\begin{figure}[!b]
\begin{center}
\setlength{\unitlength}{1cm}
\begin{picture}(8,11.6)
\put(0,8){\includegraphics[width=3.5cm]{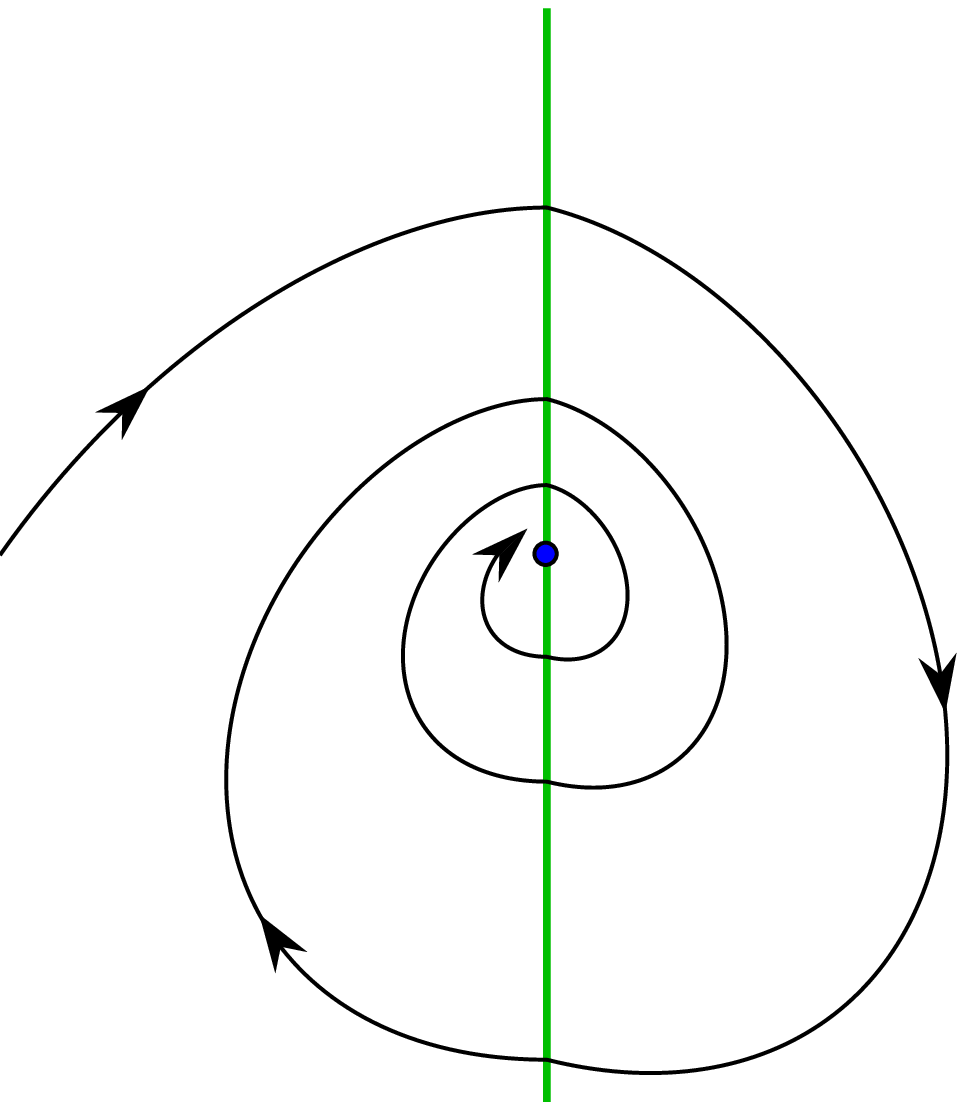}}
\put(4.5,8){\includegraphics[width=3.5cm]{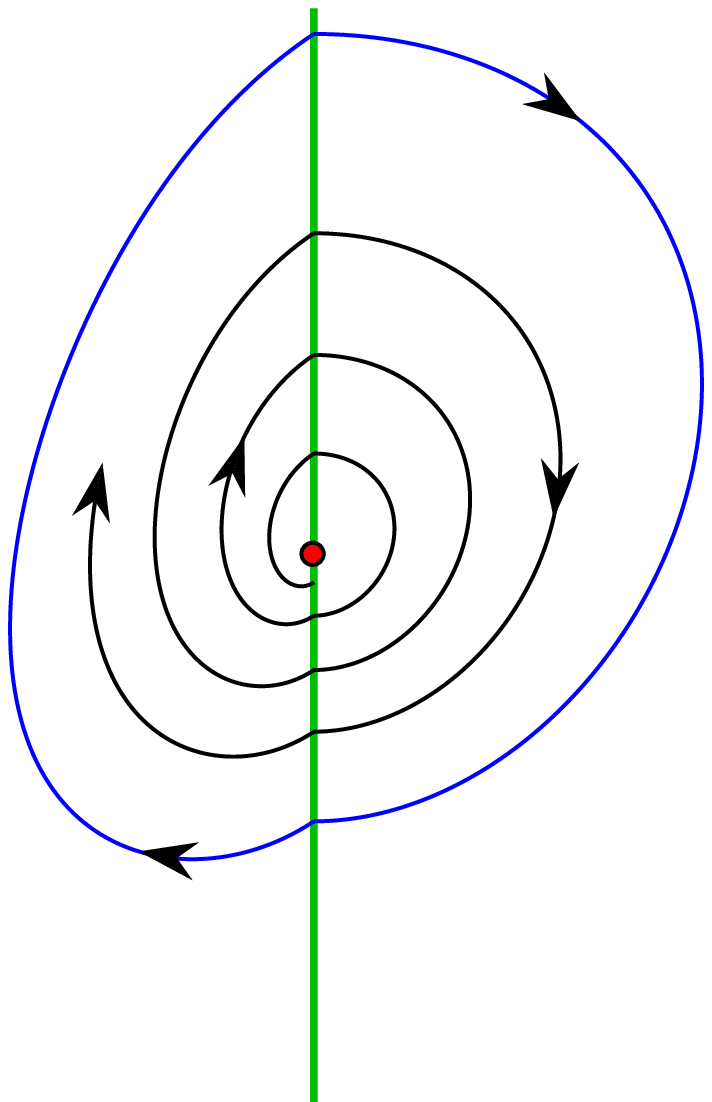}}
\put(0,4){\includegraphics[width=3.5cm]{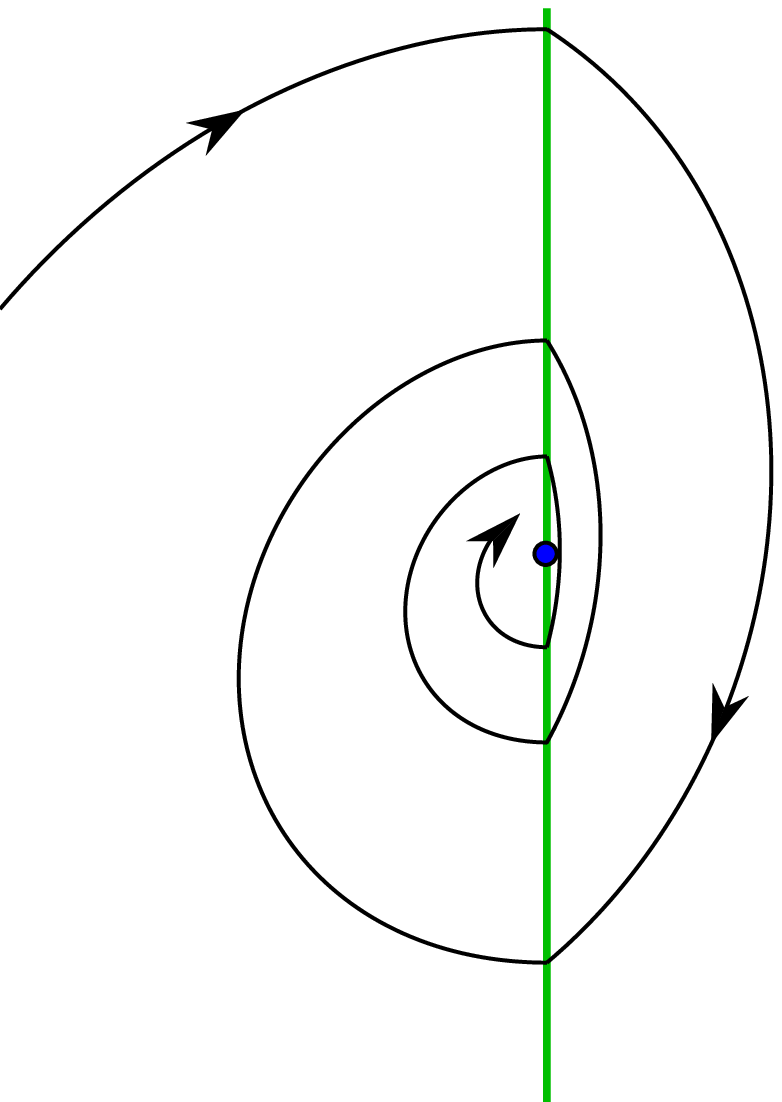}}
\put(4.5,4){\includegraphics[width=3.5cm]{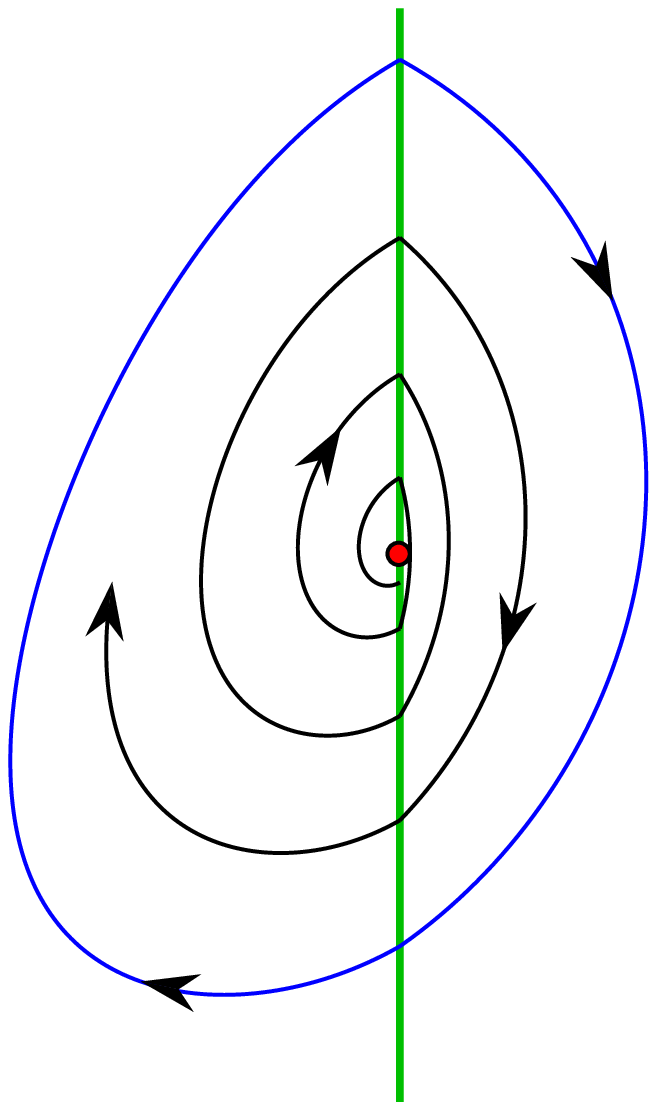}}
\put(0,0){\includegraphics[width=3.5cm]{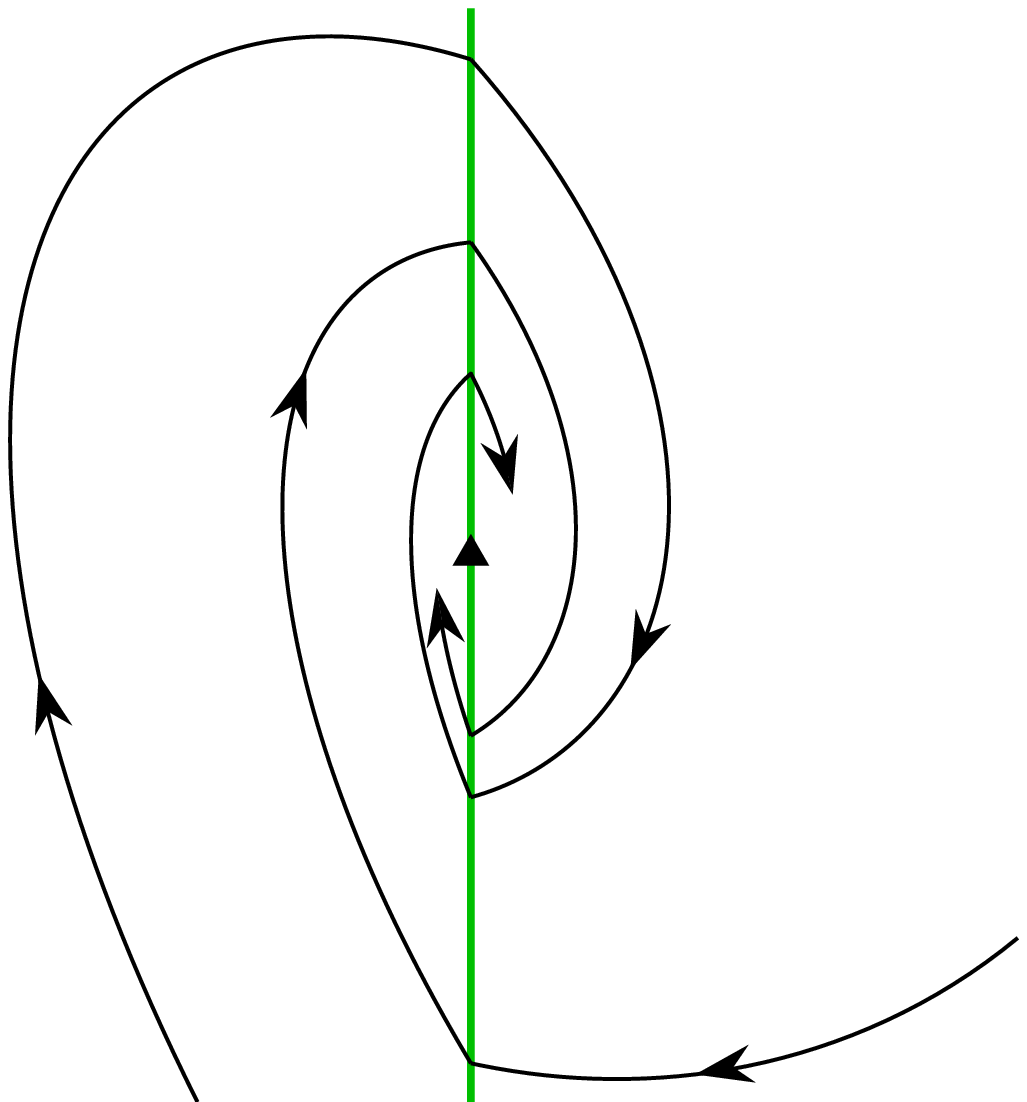}}
\put(4.5,0){\includegraphics[width=3.5cm]{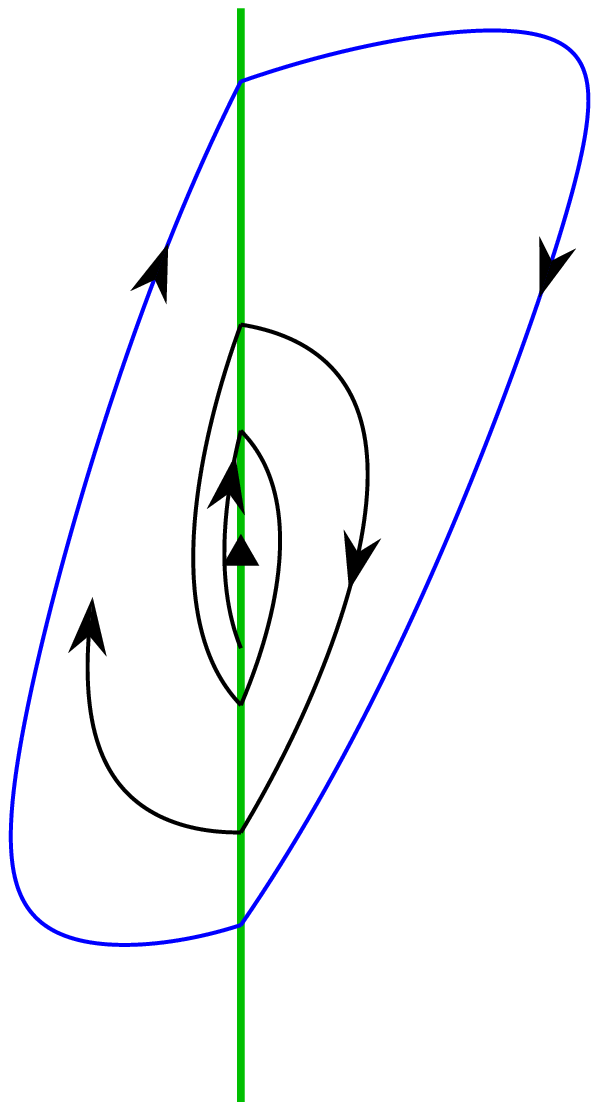}}
\put(0,11.3){\small \parbox{80mm}{\begin{center} $\circled{8}$ \end{center}}}
\put(0,7.3){\small \parbox{80mm}{\begin{center} $\circled{9}$ \end{center}}}
\put(0,3.3){\small \parbox{80mm}{\begin{center} $\circled{10}$ \end{center}}}
\end{picture}
\caption{
Phase portraits for fixed folds and foci.
\label{fig:hbl8to10}
} 
\end{center}
\end{figure}
%%%%%%%%%%%%%%%%%%%%%%%%%%%%%%%%%%%%%%%%%%%%%%%%%%%%%%%%%%%%%%%%%%%%%%%%%%%%%%%%%%%%%%%%%%%%%%%%%%%%%%%%%%%%%%%%%%%%%%%%

Now suppose $F_L$ and $F_R$ each have either a
focus or a fold fixed at the same point on $x=0$, Fig.~\ref{fig:hbl8to10}.
This point is an equilibrium, or may be treated as one,
and its stability may change, say at $\mu = 0$, as the value of $\mu$ is varied.
Generically a limit cycle is created at $\mu = 0$.
As with slipping foci and folds, there are three cases:
(i) two foci (HLB 8),
(ii) one focus and one fold (HLB 9), and
(iii) two folds (HLB 10).
% HLB 8--10.
Again only in the case of two folds does the amplitude of the limit cycle have nonlinear asymptotic growth.
The case of two foci was identified in a model of a car braking system in \cite{ZoKu06}.
In this model the amplitude of the limit cycle has square-root growth because the nonlinear terms are non-generic
(cubic instead of quadratic).

%%%%%%%%%%%%%%%%%%%%%%%%%%%%%%%%%%%%%%%%%%%%%%%%%%%%%%%%%%%%%%%%%%%%%%%%%%%%%%%%%%%%%%%%%%%%%%%%%%%%%%%%%%%%%%%%%%%%%%%%
\begin{figure}[!b]
\begin{center}
\setlength{\unitlength}{1cm}
\begin{picture}(8,15.6)
\put(0,12){\includegraphics[width=3.5cm]{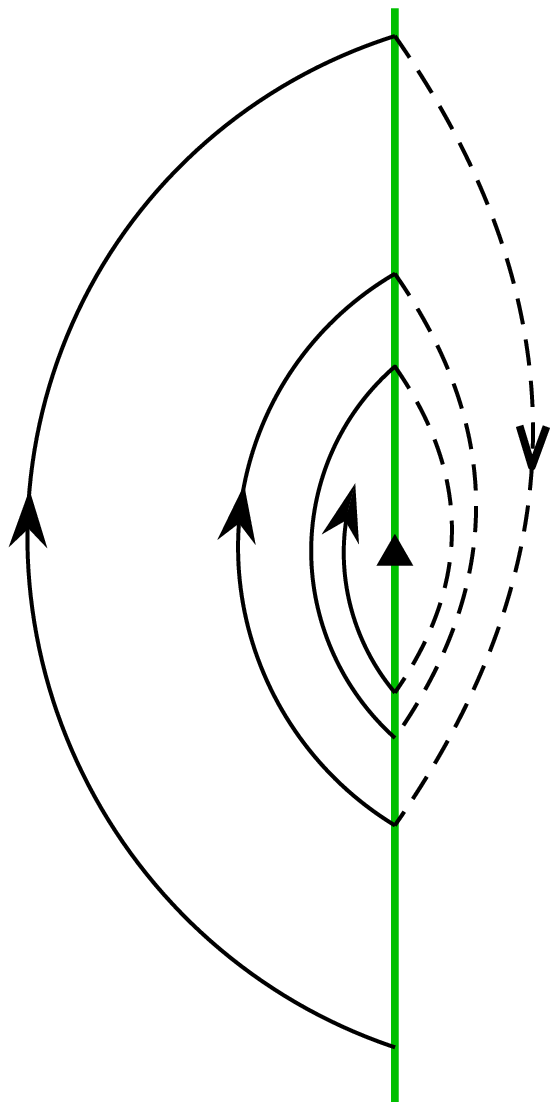}}
\put(4.5,12){\includegraphics[width=3.5cm]{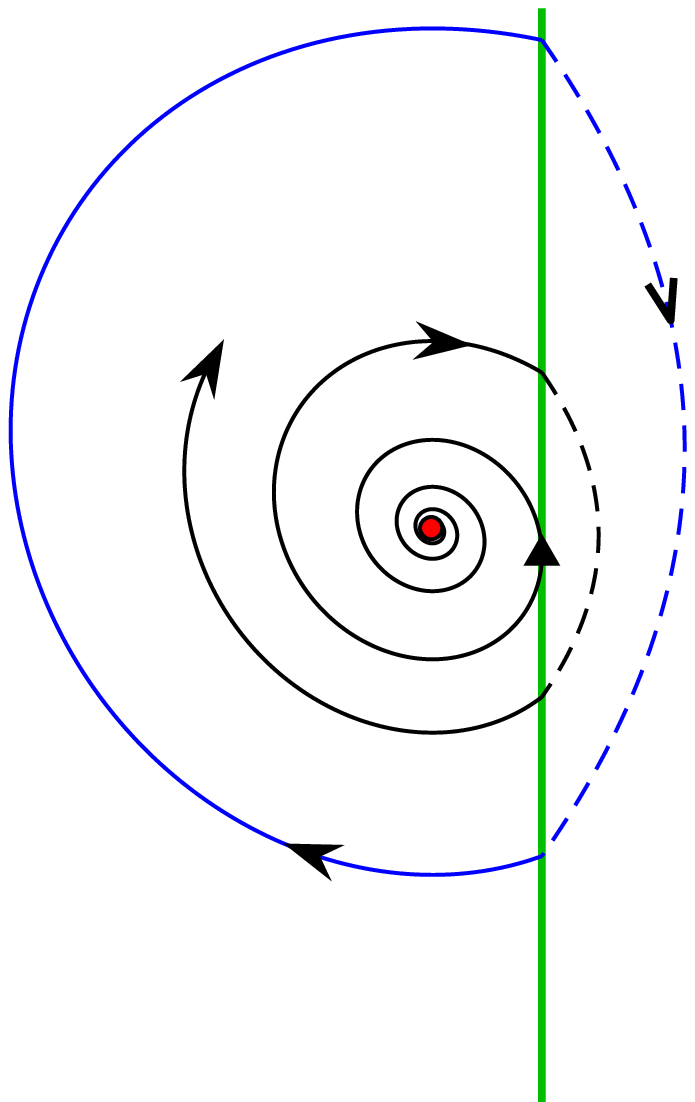}}
\put(0,8){\includegraphics[width=3.5cm]{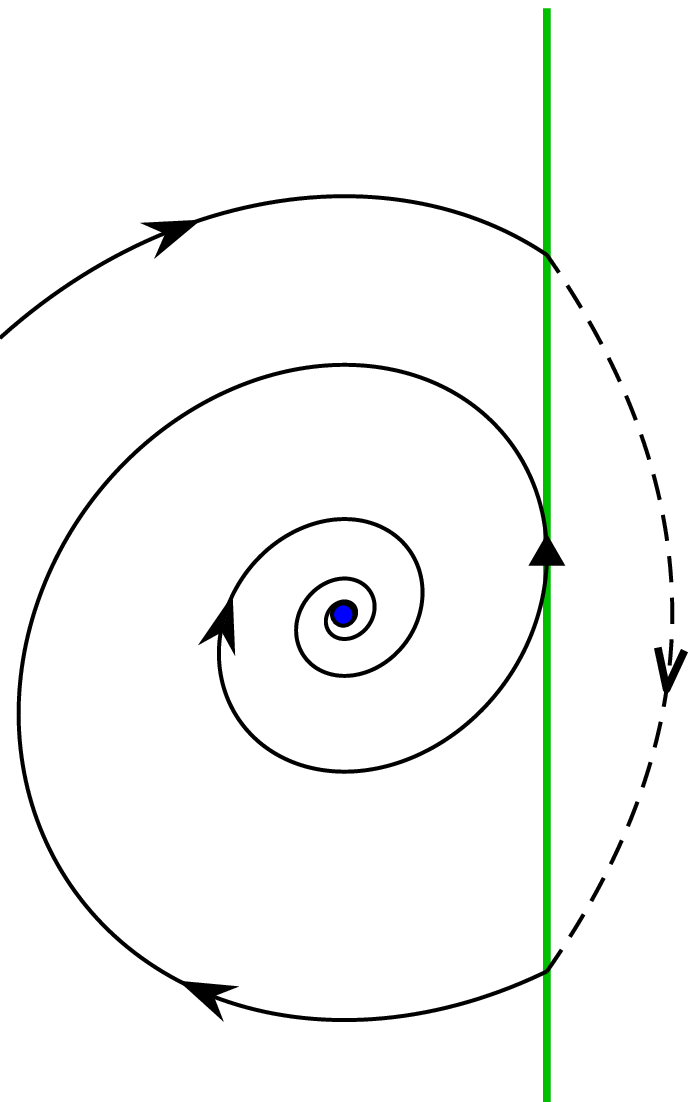}}
\put(4.5,8){\includegraphics[width=3.5cm]{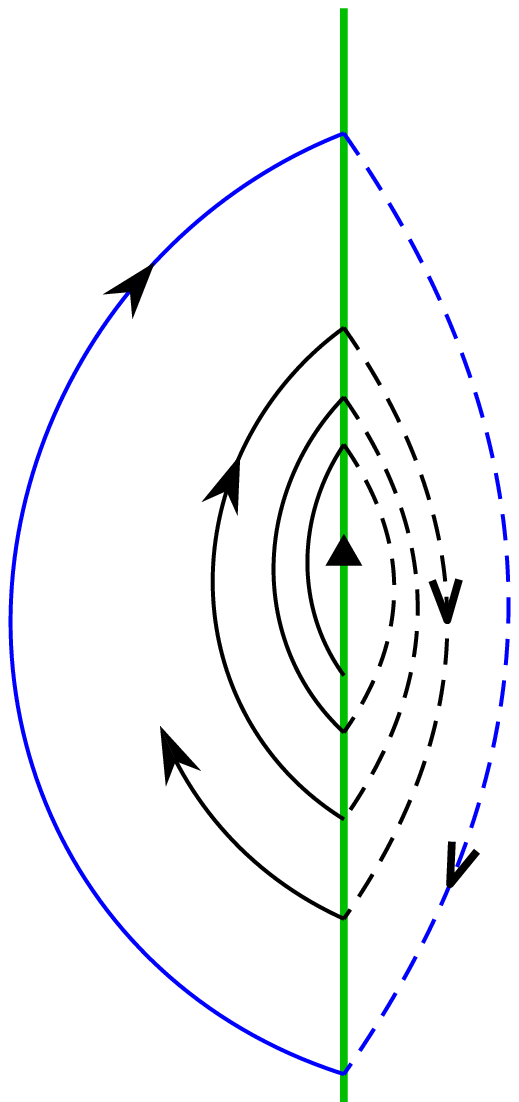}}
\put(0,4){\includegraphics[width=3.5cm]{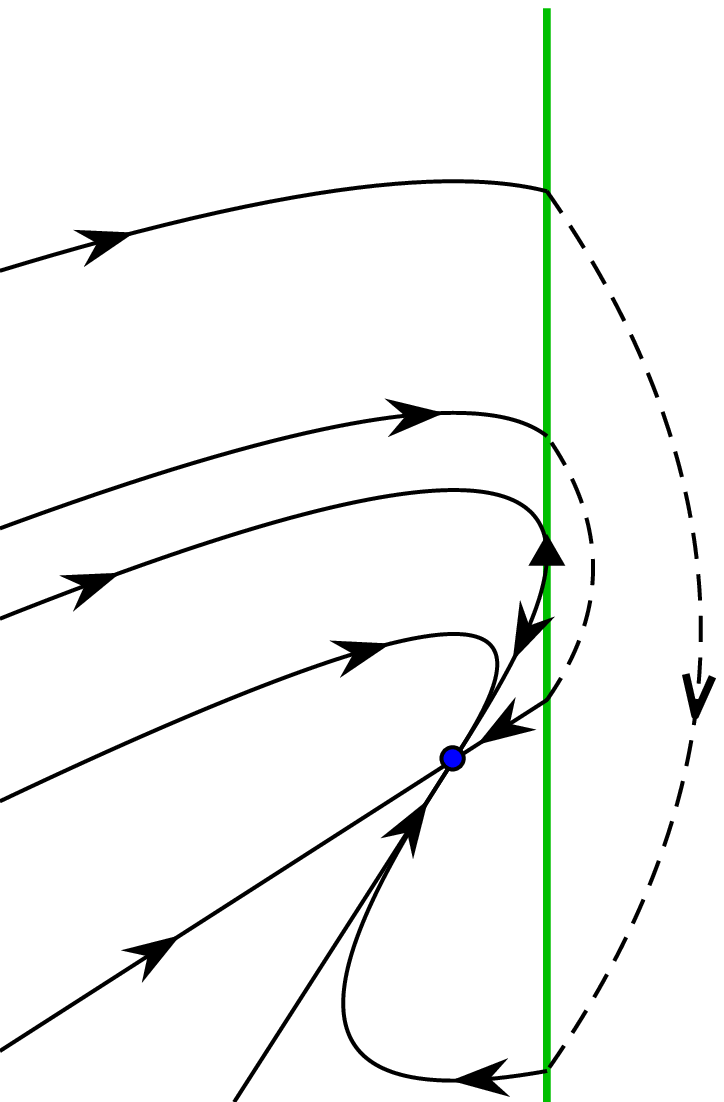}}
\put(4.5,4){\includegraphics[width=3.5cm]{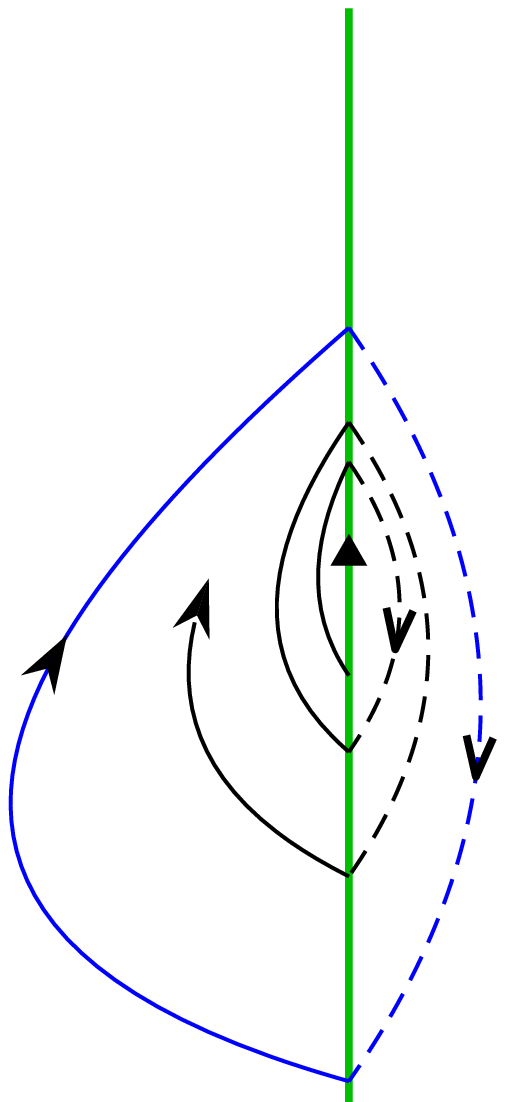}}
\put(0,0){\includegraphics[width=3.5cm]{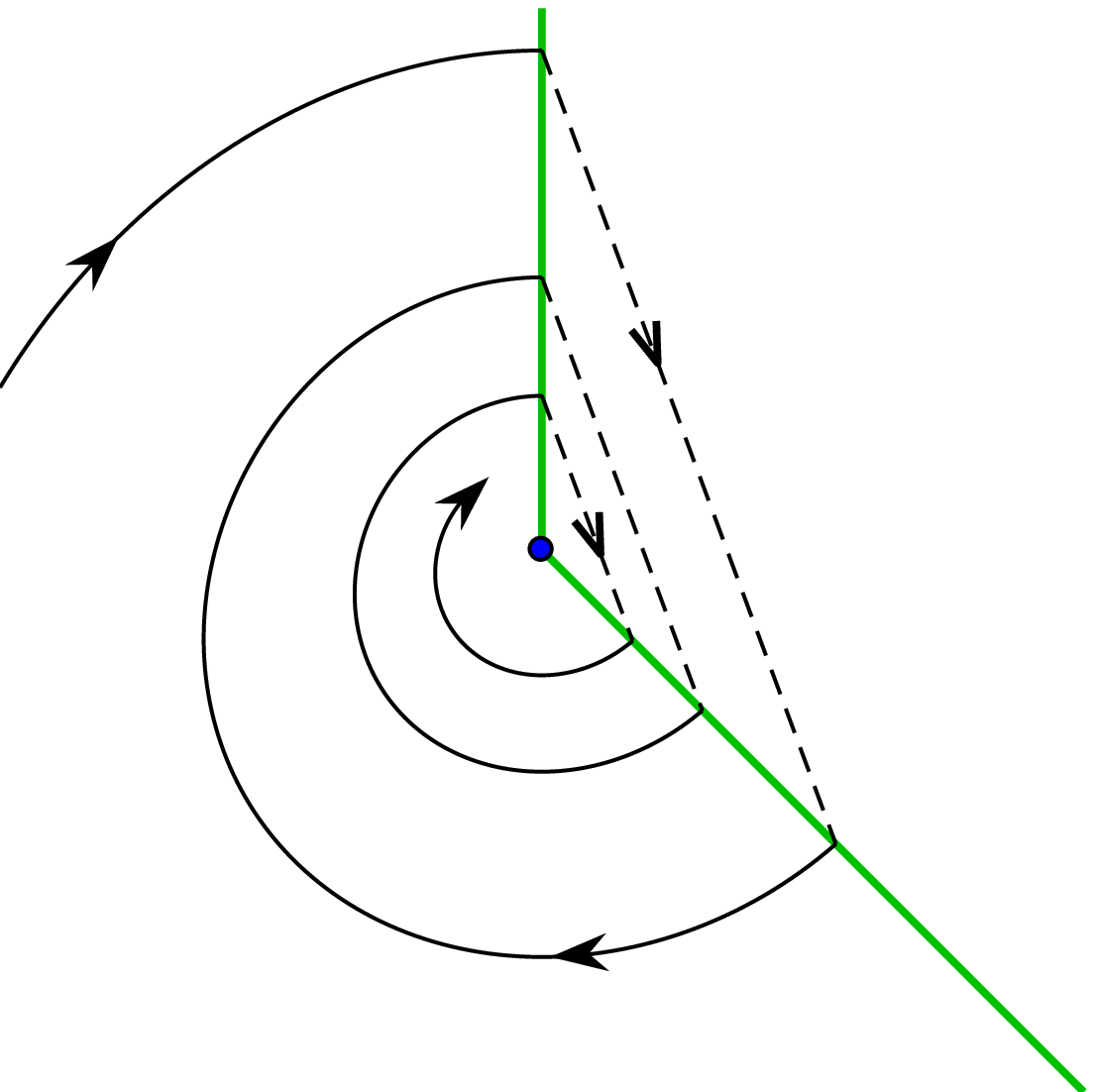}}
\put(4.5,0){\includegraphics[width=3.5cm]{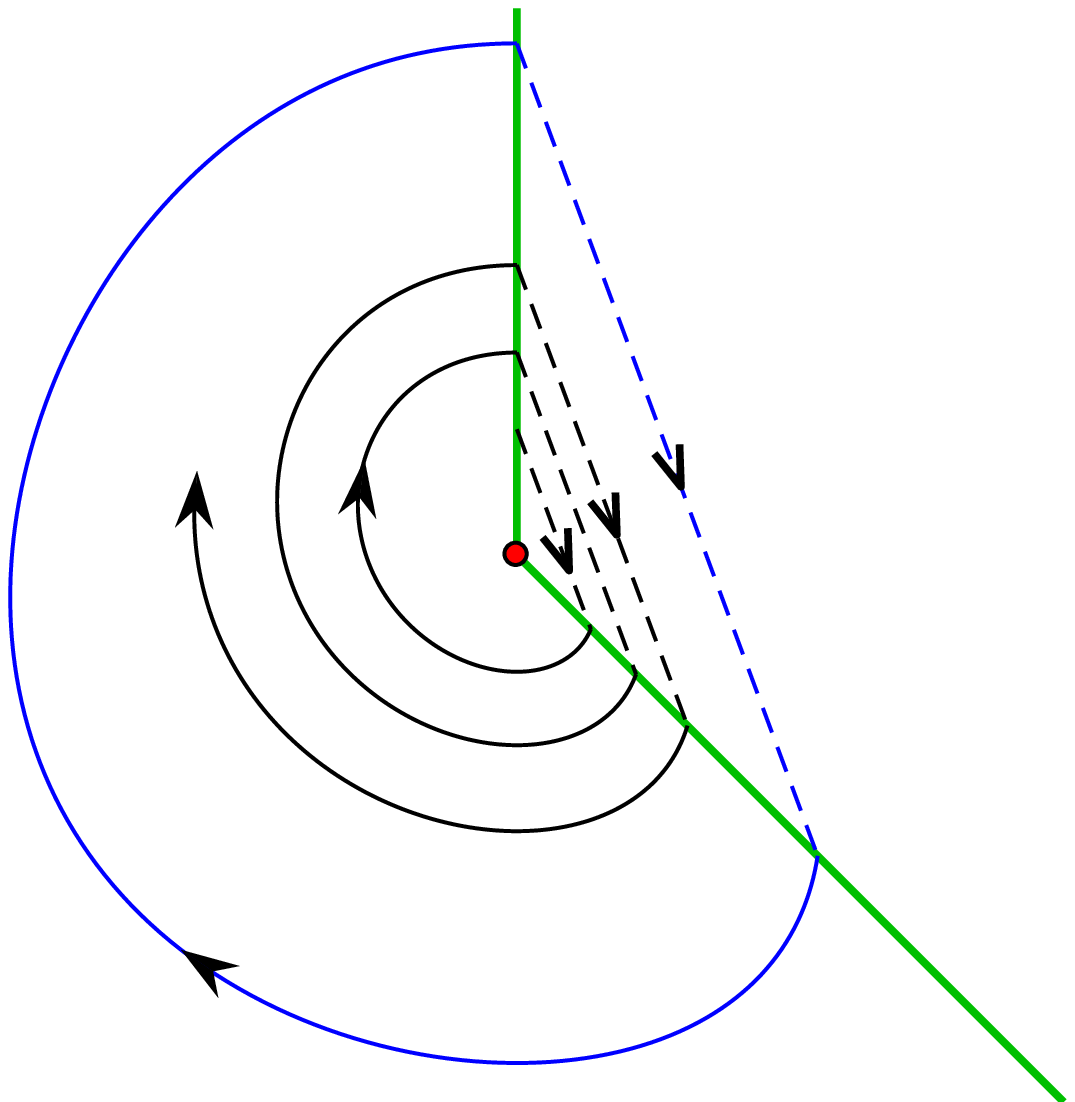}}
\put(0,15.3){\small \parbox{80mm}{\begin{center} $\circled{11}$ \end{center}}}
\put(0,11.3){\small \parbox{80mm}{\begin{center} $\circled{12}$ \end{center}}}
\put(0,7.3){\small \parbox{80mm}{\begin{center} $\circled{13}$ \end{center}}}
\put(0,3.3){\small \parbox{80mm}{\begin{center} $\circled{14}$ \end{center}}}
\end{picture}
\caption{Phase portraits for HLBs in impacting and impulsive systems.
Dashed curves indicate the action of the map $\phi$.
\label{fig:hbl11to14}
} 
\end{center}
\end{figure}
%%%%%%%%%%%%%%%%%%%%%%%%%%%%%%%%%%%%%%%%%%%%%%%%%%%%%%%%%%%%%%%%%%%%%%%%%%%%%%%%%%%%%%%%%%%%%%%%%%%%%%%%%%%%%%%%%%%%%%%%

Next we consider hybrid systems of the form
\begin{equation}
\begin{split}
\begin{bmatrix} \dot{x} \\ \dot{y} \end{bmatrix} &= F(x,y;\mu), {\rm ~for~} x < 0 , \\
y &\mapsto -\phi(y;\mu), {\rm ~when~} x = 0.
\end{split}
\label{eq:impactingODE}
\end{equation}
Let $F_1$ denote the first component of $F$.
We assume that $F_1(0,y;\mu) > 0$ for all $y > 0$, and $F_1(0,y;\mu) < 0$ for all $y < 0$.
We also assume $\phi(0;\mu) = 0$, and $\phi(y;\mu) > 0$ for all $y > 0$.
These conditions ensure that applications of the map $\phi$ are always followed by motion in $x < 0$.

Systems of this form are commonly used to model mechanical systems
with hard impacts by assuming impacting components undergo instantaneous velocity reversals \cite{WiDe00}.
In \eqref{eq:impactingODE}, $\phi$ represents the impact law and $F$ describes motion between impacts.

Suppose an equilibrium of \eqref{eq:impactingODE} collides with $x=0$
(necessarily at $x=y=0$) when $\mu = 0$, Fig.~\ref{fig:hbl11to14}.
This is a BEB and we say that the equilibrium changes from {\em admissible} (when its $x$-value is negative)
to {\em virtual} (when its $x$-value is positive).
A limit cycle is created at $\mu = 0$ in three distinct scenarios.
Specifically, a stable limit cycle can coexist with
(i) an admissible unstable focus (HLB 11),
(ii) a virtual stable focus (HLB 12),
or (iii) a virtual stable node (HLB 13).
In each case, $a = 1$ and $b = 0$, as with BEBs in continuous systems and Filippov systems.
Indeed the bifurcations can be analysed by defining a vector field in $x > 0$
that mimics the action of $\phi$ \cite{DiNo08}.

Now consider a hybrid system with a map $\phi$ from one manifold, say $x=0$, to a different manifold
(see already HLB 14 in Fig.~\ref{fig:hbl11to14}).
Systems of this form are often used to model impulsive systems,
where $\phi$ describes the action of an impulse \cite{HaCh06}.
Here we describe a HLB in such a system following \cite{Ak05}.
%Here we describe the simplest version of the HLB of \cite{Ak05}.
Suppose the system has a boundary equilibrium on $x=0$,
and that at this point the magnitude of the impulse is zero.
%at which the impulse $\phi$ is zero.
Similar to the fixed foci and folds of Fig.~\ref{fig:hbl8to10},
as parameters are varied the stability of the equilibrium can change
and a limit cycle be created (HLB 14).
A general algebraic condition determining the onset of the bifurcation is more complicated
than for HLBs 11--13, because $\phi$ provides a rotation by an arbitrary angle (not simply $180^\circ$).

%%%%%%%%%%%%%%%%%%%%%%%%%%%%%%%%%%%%%%%%%%%%%%%%%%%%%%%%%%%%%%%%%%%%%%%%%%%%%%%%%%%%%%%%%%%%%%%%%%%%%%%%%%%%%%%%%%%%%%%%
\begin{figure}[!b]
\begin{center}
\setlength{\unitlength}{1cm}
\begin{picture}(8,15.6)
\put(0,12){\includegraphics[width=3.5cm]{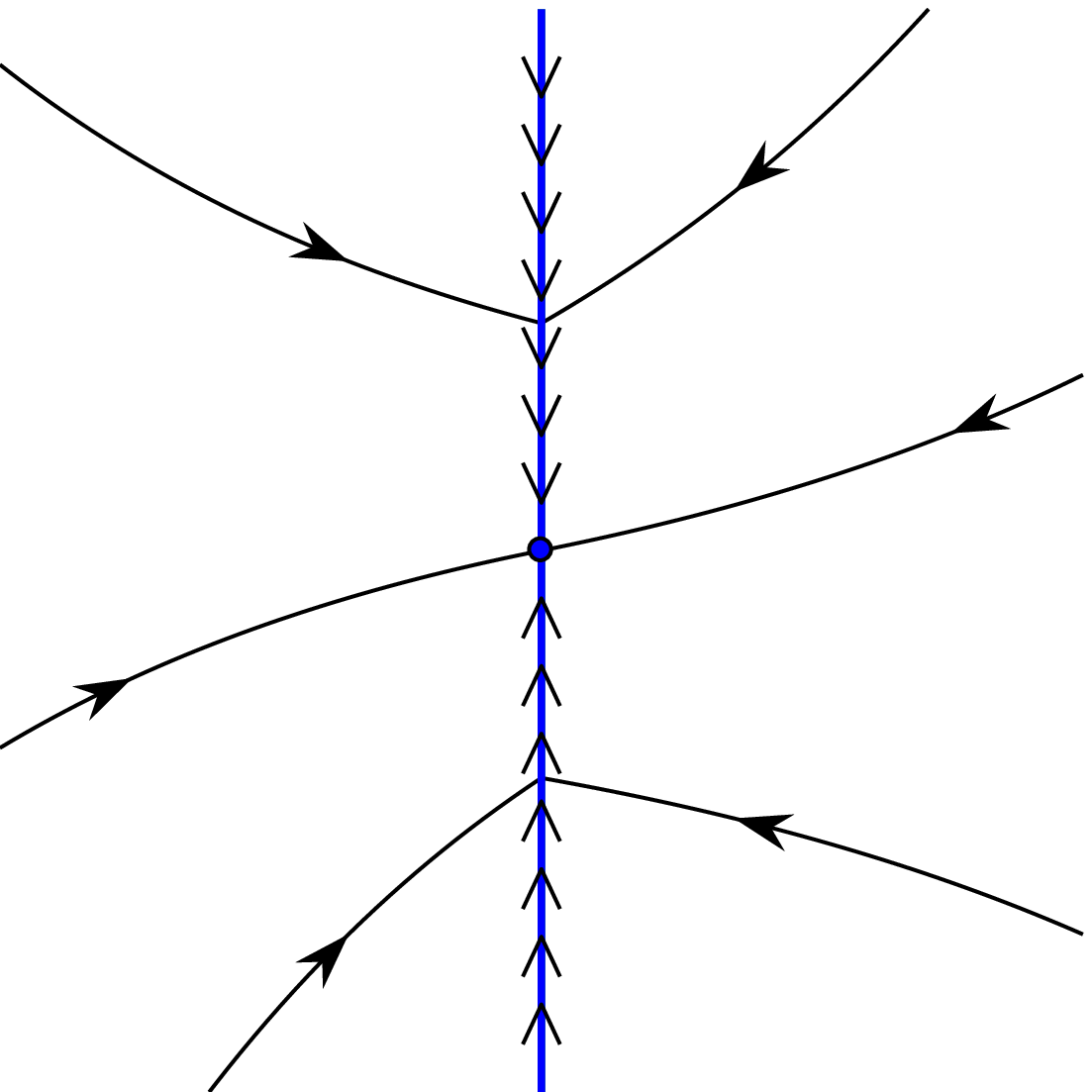}}
\put(4.5,12){\includegraphics[width=3.5cm]{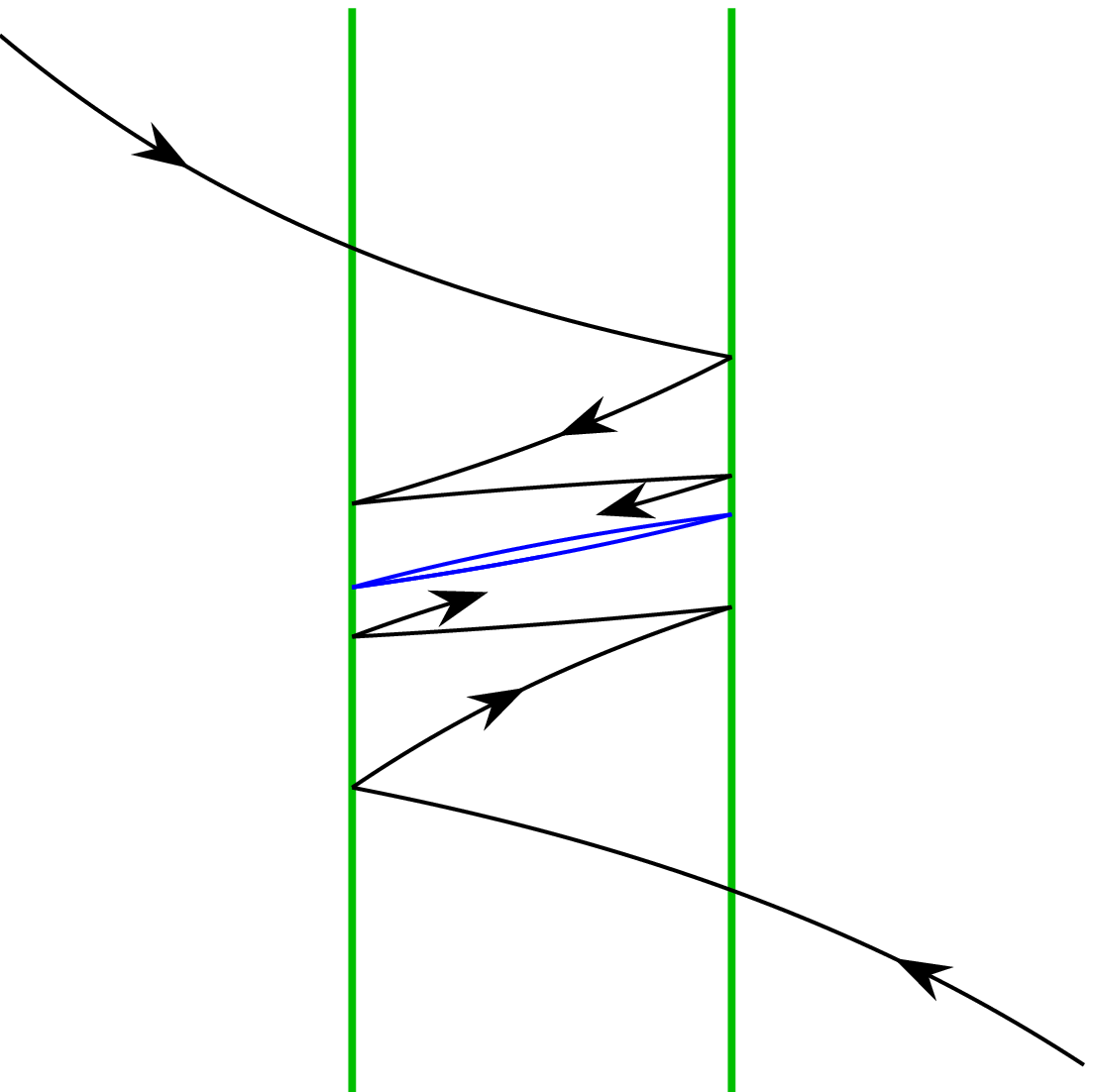}}
\put(0,8){\includegraphics[width=3.5cm]{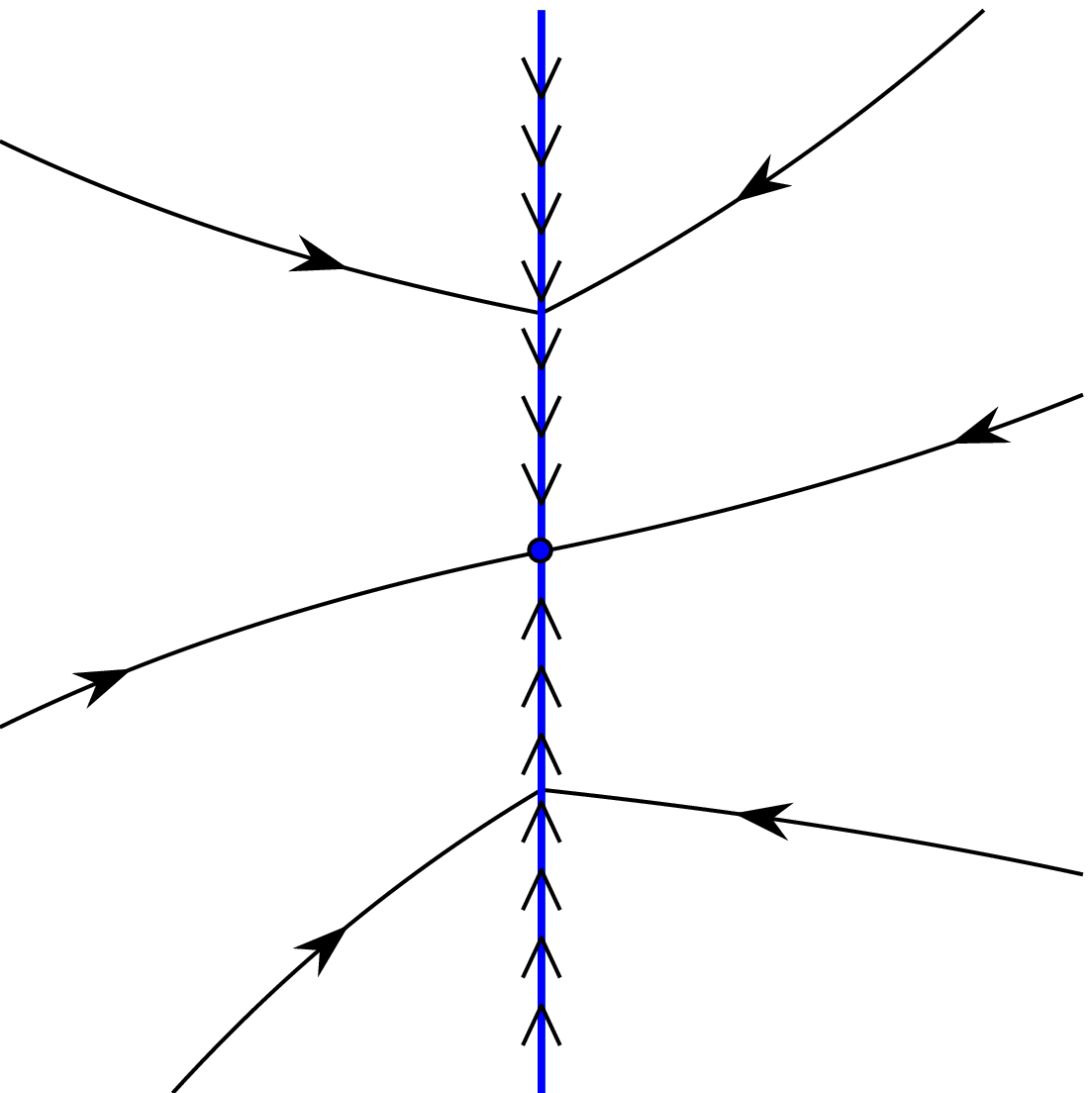}}
\put(4.5,8){\includegraphics[width=3.5cm]{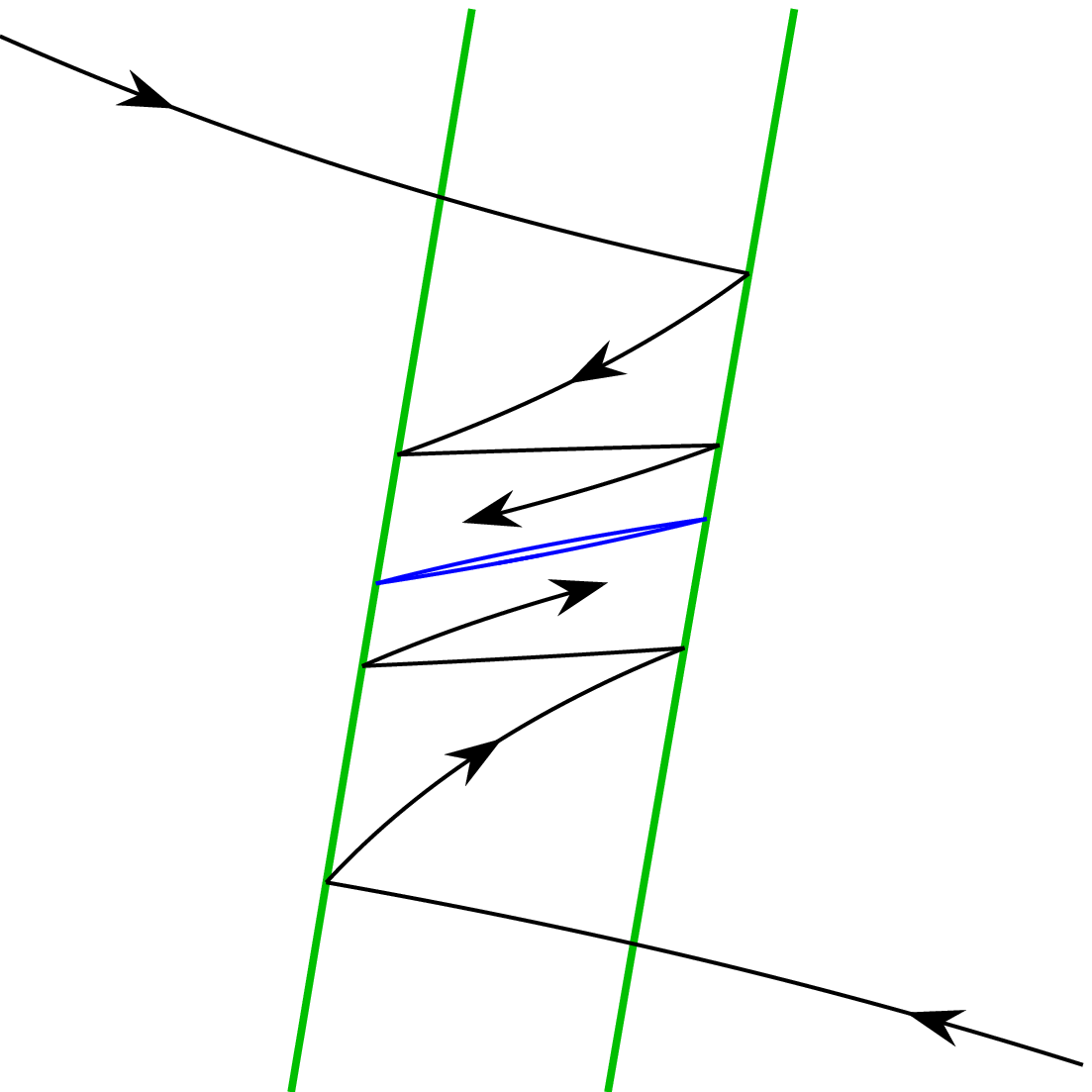}}
\put(0,4){\includegraphics[width=3.5cm]{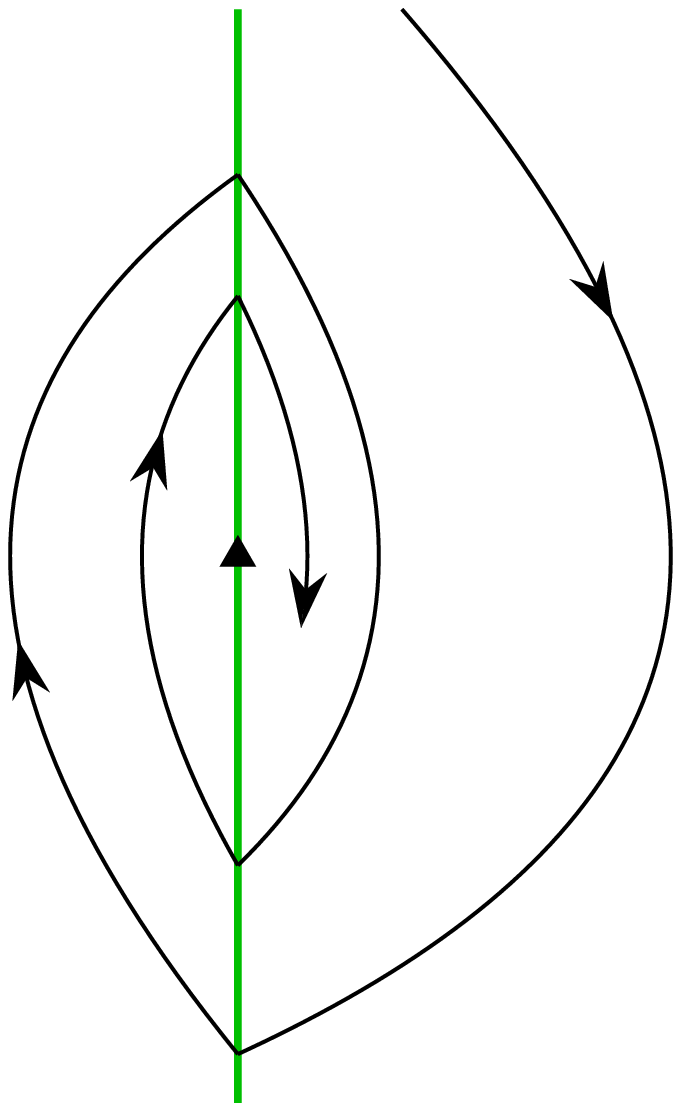}}
\put(4.5,4){\includegraphics[width=3.5cm]{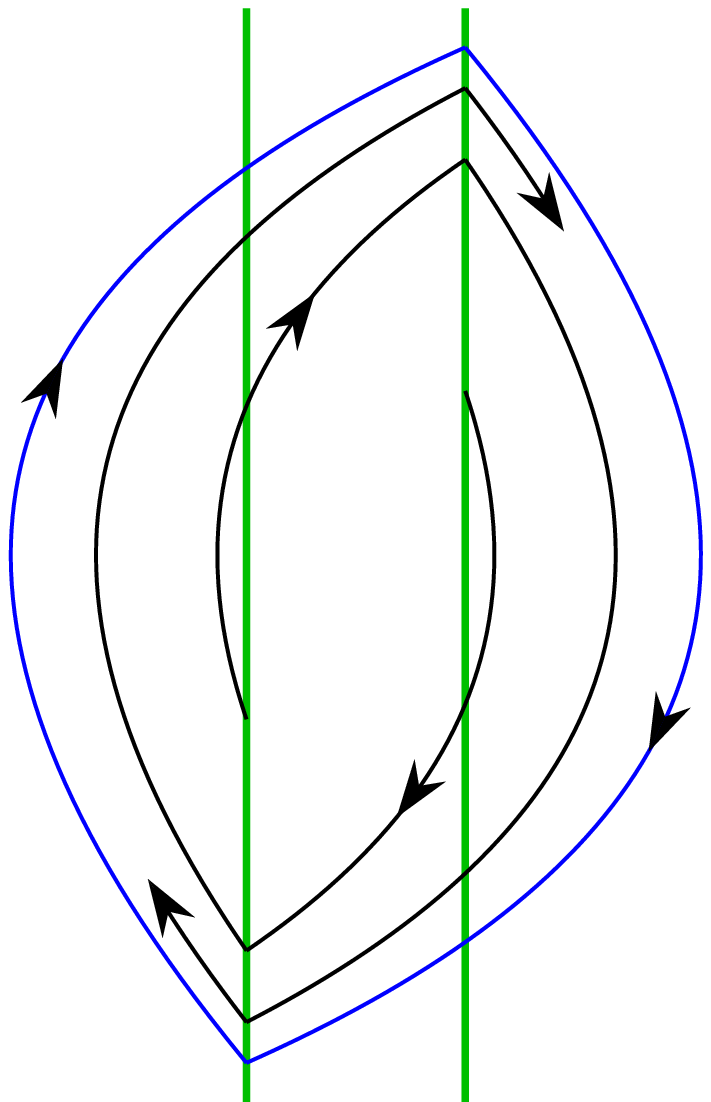}}
\put(0,0){\includegraphics[width=3.5cm]{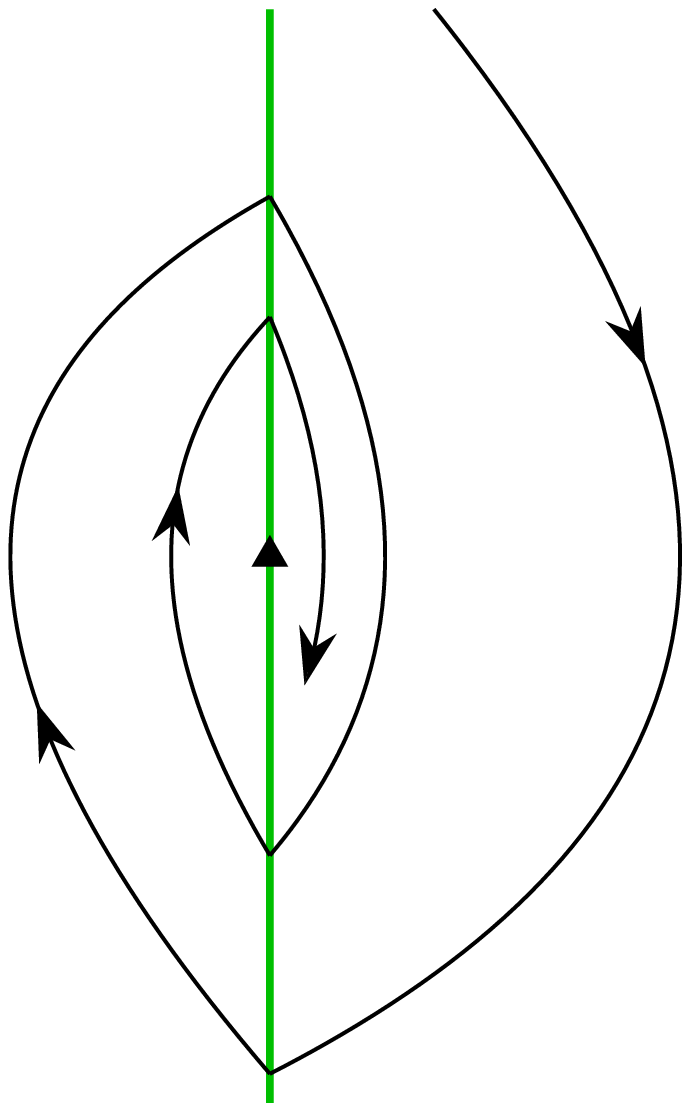}}
\put(4.5,0){\includegraphics[width=3.5cm]{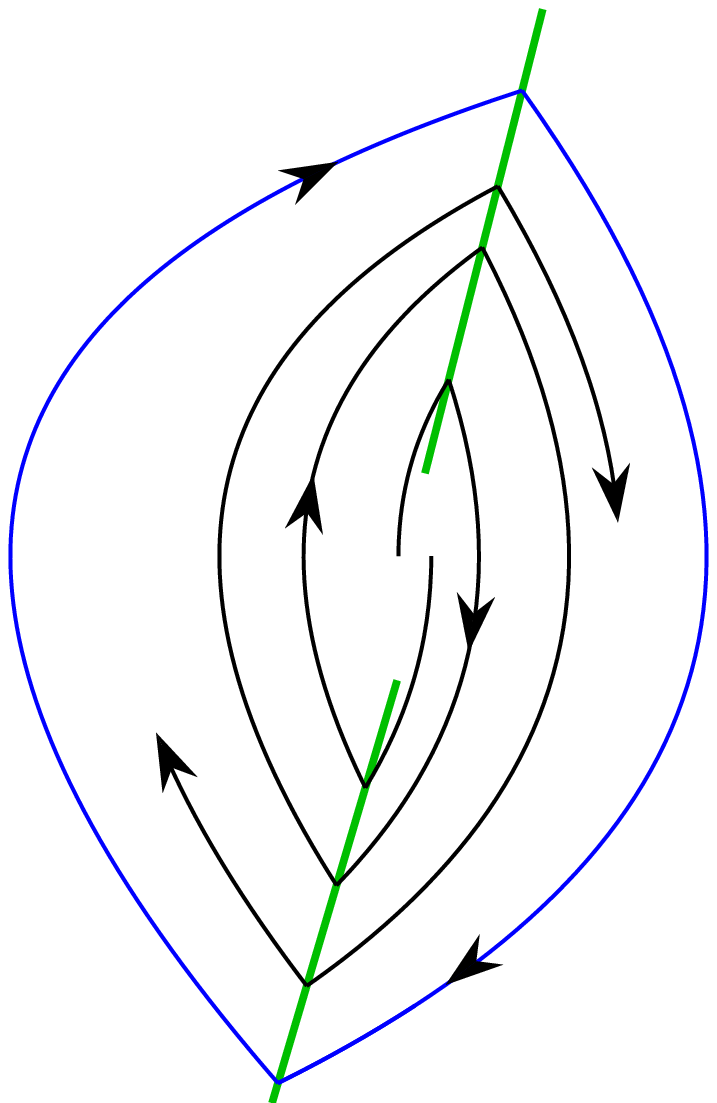}}
\put(0,15.3){\small \parbox{80mm}{\begin{center} $\circled{15}$ \end{center}}}
\put(0,11.3){\small \parbox{80mm}{\begin{center} $\circled{16}$ \end{center}}}
\put(0,7.3){\small \parbox{80mm}{\begin{center} $\circled{17}$ \end{center}}}
\put(0,3.3){\small \parbox{80mm}{\begin{center} $\circled{18}$ \end{center}}}
\end{picture}
\caption{Phase portraits for systems with hysteresis or time-delay.
\label{fig:hbl15to18}
} 
\end{center}
\end{figure}
%%%%%%%%%%%%%%%%%%%%%%%%%%%%%%%%%%%%%%%%%%%%%%%%%%%%%%%%%%%%%%%%%%%%%%%%%%%%%%%%%%%%%%%%%%%%%%%%%%%%%%%%%%%%%%%%%%%%%%%%

%Next we consider perturbations of Filippov systems.
Next we introduce perturbations.
Filippov systems are useful mathematical models of switched systems,
particularly control systems and electrical systems, when the time between switches
is small relative to the overall time-scale of the dynamics.
Such models may be made more realistic by incorportating
hysteresis or time-delay to capture individual switching events.
This regularises the switching manifold, replacing sliding motion with rapid switching.

If a Filippov system of the form \eqref{eq:FilippovBEB}
has a stable pseudo-equilibrium at $(x,y) = (0,0)$, say,
then upon replacing the switching condition at $x = 0$ with hysteretic conditions at $x = \pm \mu$,
we generate a limit cycle \cite{MaHa17} (HLB 15).
By instead supposing that orbits switch at a time $\mu > 0$ after crossing $x=0$,
we also generate a limit cycle (HLB 16).
In both cases, the period of the limit cycle is asymptotically proportional to $\mu$ (i.e.~$b = 1$).

Now suppose \eqref{eq:FilippovBEB} has a stable invisible-invisible two-fold at $(x,y) = (0,0)$,
see Fig.~\ref{fig:hbl15to18}.
By adding hysteresis or time-delay as above we again generate a stable limit cycle.
Interestingly, hysteresis (HLB 17) gives $a = b = \frac{1}{3}$ \cite{Ma17}.
This is because the Poincar\'e map on $x = -\mu$ has the form
\begin{equation}
P(y;\mu) = \sqrt{y^2 + c_1 \mu + c_2 y^3 + \cdots} \,,
\nonumber
\end{equation}
and so the solution to the fixed point equation, $y = P(y;\mu)$, involves a cube-root.
%$y = \left( \frac{-c_1 \mu}{c_2} \right)^{\frac{1}{3}} + \cdots$.
Time-delay (HLB 18) instead gives $a = b = \frac{1}{2}$,
as with earlier HLBs involving two folds.

%%%%%%%%%%%%%%%%%%%%%%%%%%%%%%%%%%%%%%%%%%%%%%%%%%%%%%%%%%%%%%%%%%%%%%%%%%%%%%%%%%%%%%%%%%%%%%%%%%%%%%%%%%%%%%%%%%%%%%%%
\begin{figure}[!h]
\begin{center}
\setlength{\unitlength}{1cm}
\begin{picture}(8,7.6)
\put(0,4){\includegraphics[width=3.5cm]{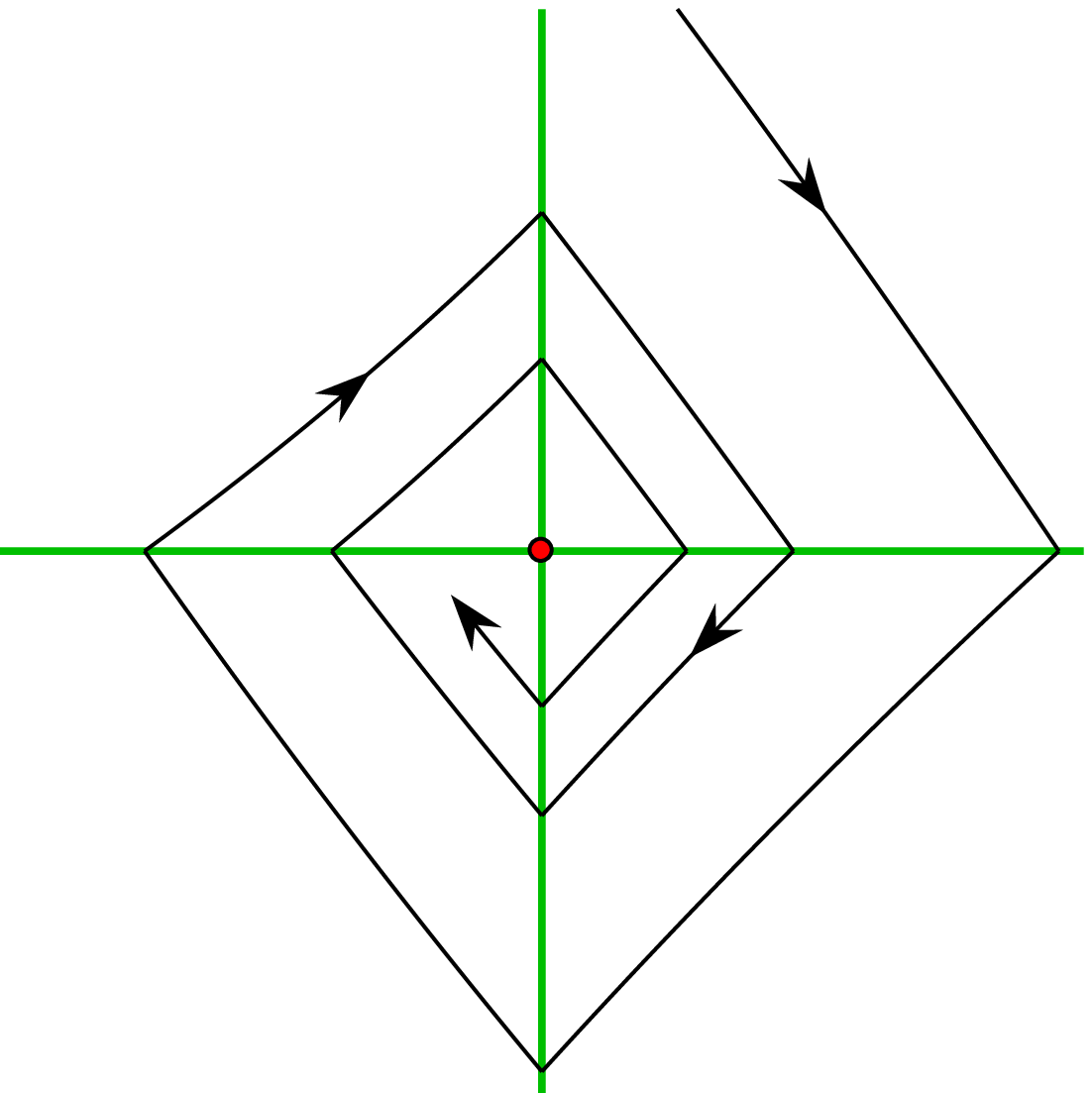}}
\put(4.5,4){\includegraphics[width=3.5cm]{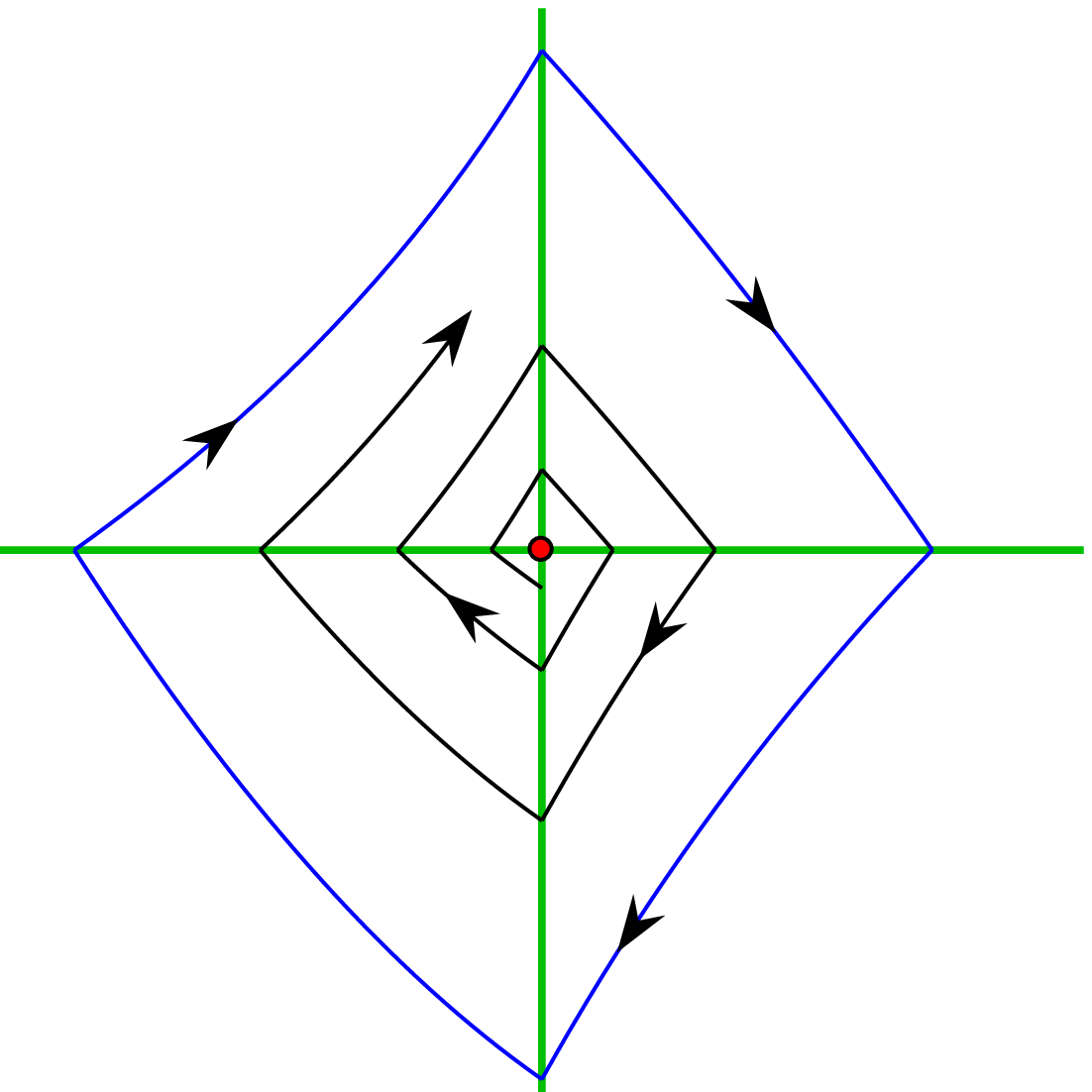}}
\put(0,0){\includegraphics[width=3.5cm]{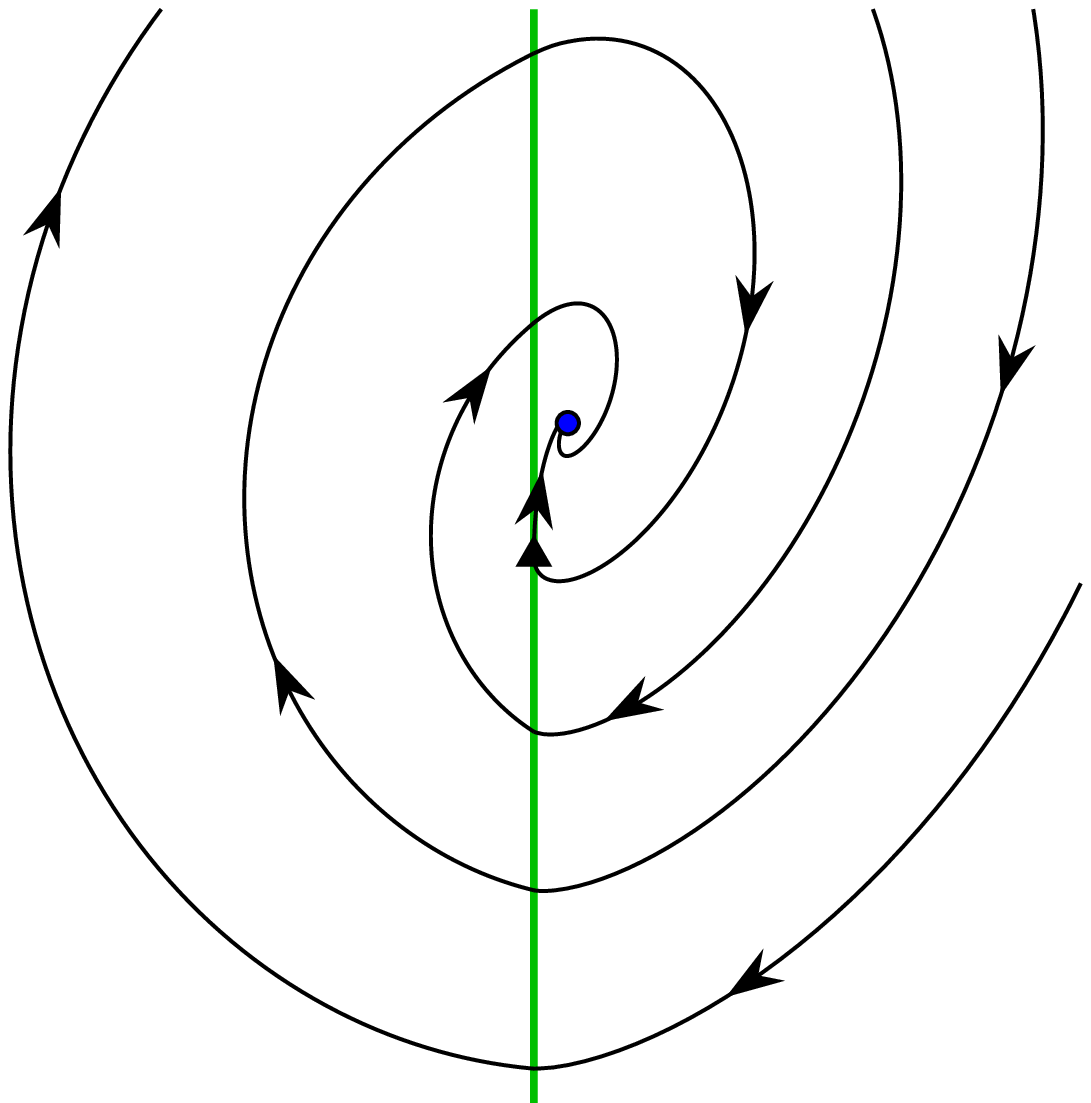}}
\put(4.5,0){\includegraphics[width=3.5cm]{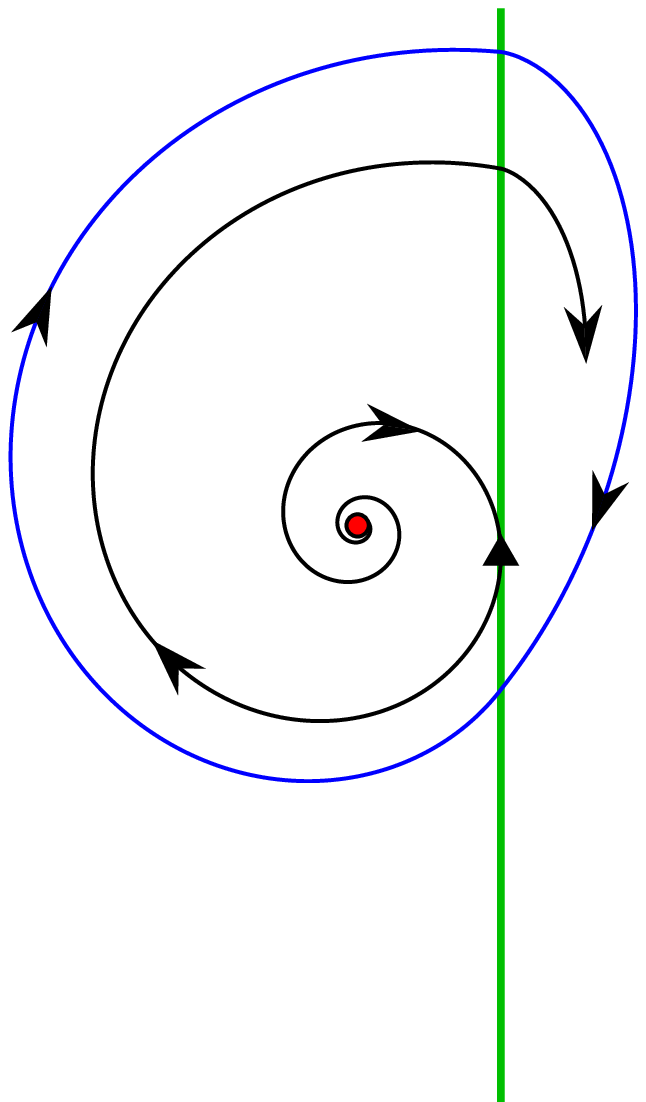}}
\put(0,7.3){\small \parbox{80mm}{\begin{center} $\circled{19}$ \end{center}}}
\put(0,3.3){\small \parbox{80mm}{\begin{center} $\circled{20}$ \end{center}}}
\end{picture}
\caption{
%Phase portraits for a system with intersecting switching manifolds, HLB 19,
%and a system with a square-root singularity, HLB 20.
Phase portraits for systems with intersecting switching manifolds
and a square-root singularity.
\label{fig:hbl19to20}
} 
\end{center}
\end{figure}
%%%%%%%%%%%%%%%%%%%%%%%%%%%%%%%%%%%%%%%%%%%%%%%%%%%%%%%%%%%%%%%%%%%%%%%%%%%%%%%%%%%%%%%%%%%%%%%%%%%%%%%%%%%%%%%%%%%%%%%%

Now we consider a Filippov system with two switching manifolds
(modelling, say, a switched system with two independent switches).
%Switched systems with multiple independent switches are naturally modelled
%by piecewise-smooth systems with multiple switching manifolds.
Since switching manifolds are codimension-one surfaces,
in two dimensions two switching manifolds generically intersect at a point.
This point may behave like an equilibrium, and, as with HLBs 8--10,
under parameter variation the intersection point cay change stability and a limit cycle be created.
This is HLB 19, see Fig.~\ref{fig:hbl19to20}.
This requires orbits to spiral around the intersection point and occurs in the neuron model of \cite{HaEr15},
where two different discontinuous functions are used to model the firing rates of excitatory and inhibitory neurons.
%The bifurcation requires orbits to spiral around the intersection point,
%and stability can be determined by computing the motion of a full revolution 

Finally, in \cite{NiCa16} the authors study BEBs in a piecewise-smooth continuous neuron model
that involves a square-root singularity.
By this we mean that, in the form \eqref{eq:pwscODE}, one component, say $F_R$, has a $\sqrt{|x|}$-term.
As in the usual continuous scenario (HLBs 1--2),
both $F_L$ and $F_R$ have an equilibrium, and these coincide at the BEB.
Here, however, in order to generate a limit cycle locally, both equilibria must be foci (HLB 20).
The bifurcation therefore closely resembles HLB 1,
however the square-root term prevents the limit cycle from deeply penetrating the $x > 0$ half plane.
While the amplitude of the limit cycle is asymptotically proportional to $\mu$,
its maximum $x$-value is asymptotically proportional to $\mu^2$.
For this reason, HLB 20 is perhaps best viewed an intermediary of HLB 1 and HLB 3.

The HLBs presented here are not intended to form a complete list
but hopefully cover the most fundamental scenarios and those reported in mathematical models.
Other bifurcations of piecewise-smooth systems that involve limit cycles,
but less closely resemble a Hopf bifurcation, include discontinuity induced bifurcations
at which two limit cycles are created simultaneously \cite{MoBu14,SiKu12b}.
BEBs can mimic saddle-node bifurcations
in that two equilibria (one of which is a pseudo equilibrium in the case of Filippov systems)
may collide and annihilate at the bifurcation.
Here a local limit cycle is created at the same time
if the limiting piecewise-linear system satisfies a certain global property \cite{SiMe12,KuRi03}.
Limit cycles can also be created in global bifurcations % unique to piecewise-smooth systems,
such as  `canard super-explosions' \cite{DeFr13,RoCo12},
and certain bifurcations of piecewise-linear systems with three or more components \cite{LlPo15,PoRo15}.

%{\footnotesize
%\bibliographystyle{unsrt}
%\bibliography{../DynSyst,../MathBio,../Misc,../OtherTheory,../PWS,../Stoch}
%}

\end{document}